\documentclass[11pt]{article}
\normalfont
\renewcommand{\baselinestretch}{0.98}

\usepackage{amsfonts}
\usepackage{times}

\newcounter{myfn}[page]
\renewcommand{\thefootnote}{\fnsymbol{footnote}}
\newcommand{\myfootnote}[1]{\setcounter{footnote}{\value{myfn}}%
    \footnote{#1}\stepcounter{myfn}}

\newcommand{\myDATE}{November 3, 2008}
\newcommand{\ClarifyingRemark}{(See \S 4 for clarifying comments on the 
results presented here.)}

\newcommand{\statset}[3]{\hspace*{0.1in}\parbox{0.4in}{\hfill 
\footnotesize #1} \ 
\parbox{0.4in}{\hfill \footnotesize #2}\ \parbox{0.4in}{\hfill 
\footnotesize #3}}

\begin{document}

\thispagestyle{empty}

\vspace*{-0.7in}
\hfill {\footnotesize November 20, 2008}

\begin{center}
{\large \bf A view of 
mathematics research productivity at U.S.\ regional public universities}

\renewcommand{\thefootnote}{1}
Robert G.\ Donnelly\footnote{Email: 
{\tt rob.donnelly@murraystate.edu},\ Fax: 
1-270-809-2314}

\vspace*{-0.075in}
Department of Mathematics and Statistics, Murray State
University, Murray, KY 42071 USA
\end{center}

\vspace*{-0.25in}
\begin{abstract}
Statistical summaries of certain kinds of 
mathematics research output are given for a large sample of U.S.\ 
regional public universities.  
These statistical summaries are reported using a variety of 
metrics that distinguish between single-authored and collaborative 
work and account for publication length.   
\end{abstract}

\vspace*{-0.1in}
{\bf \S 1. Introduction.}  
In June 2006, the 
Minister of State of Higher Education in Britain proposed 
a regimen for reforming the system of assessment and evaluation 
of research productivity in higher education \cite{ReformProposal}.  
The new system would 
involve automated assessment of data.  The proposal stated 
that the old system 

\vspace*{0.1in}
\hspace*{0.5in}\parbox{5.5in}{\footnotesize 
should be replaced with a new and lighter-touch system based largely 
on metrics.  The principle of using information that is already 
collected routinely to assess research quality and allocate funding 
must be the right one.  The savings of time and effort that this can 
bring for university teachers and administrators alike should be 
welcomed by all, as should the transparency that a system based on 
publicly available data potentially offers.}    

\vspace*{0.1in}
In October 2006, the Council for the Mathematical Sciences (or CMS, comprised 
of the Institute for Mathematics and its Applications, the London 
Mathematical Society, and the Royal Statistical Society) responded to 
the reform proposal \cite{CMSResponse}, \cite{AppendixB}.  
While acknowledging that the use of a variety of 
numerical measures is an important 
part of the picture of productivity, 
the CMS expressed serious concerns that excessive reliance on numbers can 
adversely affect the way mathematics research is practiced and that 
metrics that work well for other science disciplines might not work 
well for mathematics.  
In December 2006, it was announced that mathematics would, for the 
time being, be excluded 
from the new ``lighter-touch'' system in favor of a 
regimen that would rely both on metrics and peer assessment 
\cite{HEFCEAnnouncement}. 

These events 
partly 
inspired the study presented in this paper 
of mathematics research output at 
U.S.\ regional public universities.  This study was undertaken with the 
following goals: (1) To provide some interesting and hopefully useful 
summary data on mathematics research output as evidenced by 
publications; 
(2) to support the use of a variety of measures in   
assessments of mathematics faculty research output rather 
than the flawed metric that is often exclusively used (that 
is, length of the publication list); (3) to demonstrate that raw data 
supporting such measures is widely available, although some effort 
might be required to collect and process the data; and (4) to promote 
the view that in any assessment of mathematics faculty performance, numerical 
measures should be used with great care and 
qualified by the perspectives of knowledgeable and established 
practitioners.    

This fourth point is perhaps most important.  The following concerns 
should be kept in mind.  
The existence of such metrics should not 
diminish the value or discourage the pursuit 
of other kinds of 
scholarly mathematical activities which are beyond the scope of this 
study (such as time and energy given to teaching and working with students, 
course development, 
direction of theses and other student projects, 
talks and presentations, 
reviewing and refereeing, 
editorial work, 
expository writing, 
notes on pedagogy, 
textbook authoring, 
conference organizing, 
service to scholarly organizations, 
consulting, translating, grant writing, etc).  
Moreover, it should be noted that 
it is difficult for statistics about output volume 
to account for less tangible 
factors such as quality. 
These concerns are easily forgotten  
when undue prominence is given to numerical measures, and this can 
consequently motivate 
behavior that undermines the efficacy of the discipline.  
The CMS voiced similar concerns in 
\cite{CMSResponse}, advocating that metrics for assessing research 
productivity should not 
be considered in isolation without the context and moderation of 
expert viewpoints. 
Nonetheless, as the CMS readily acknowledges \cite{CMSResponse}, 
\cite{AppendixB}, 
a variety of reliable metrics is desirable. 

There is now an extensive literature on what is known as 
``evaluative bibliometrics.'' 
Since the 1960's, 
databases of publication information have been used to support metrics 
for assessments of research performance.  
Such metrics typically measure output volume or the 
impact of research and are applied variously to individuals, 
programs, journals, disciplines, countries, or regions of the world.  
At this time it appears 
there is no uniform standard for measuring scientific 
research output, see 
\cite{Larsen}.  For one possible approach to the problem of 
standardizing counting methods, see \cite{Gauffriau}.  
One critical issue is whether ``total counting'' or ``fractional 
counting'' should be used in assessing numbers of papers: roughly 
speaking, total counts do not distinguish 
between single-authored and collaborative 
publications, whereas fractional counts do. 
Impact is often measured by citations. 
The so-called ``Hirsch index'' 
(or $h$-index) introduced by Hirsch \cite{Hirsch} combines output volume and 
impact in a single numerical measure.  
This metric and its variations are now widely used but have their 
limitations, as 
acknowledged by Hirsch and pointed out by others (see for example 
\cite{Bornmann}).  
Many mathematical societies have pointed out the limited usefulness 
of citation statistics (such as 
the Thomson Scientific impact factors) for measuring the impact of 
mathematics research (\cite{IMU}, \cite{BMS}, \cite{AppendixB}).  
Partly this is due 
to different citation habits in mathematics,  
to overall publication rates, and to the nature of the discipline 
itself.  
As noted in  
\cite{Iglesias} and \cite{IglesiasArxiv}, it does seem that 
``\ldots`Mathematics' is a 
field that has quite specific rules, and probably requires 
individualized treatment.''  
Since the main interest of this study is  
output volume of mathematics research, 
measures of impact will not be considered here.  

The remainder of the paper aims to address goals (1), (2), and (3). 
In \S 2, the design and methods of the study are described.  
A key aspect of the study is noted in this section, namely 
the high accessibility of the web-based resources that provided 
the input data for the study. 
In \S 3, descriptions are given for 
the metrics used here.  
In measuring output volume, these metrics notably 
distinguish between 
single-authored and collaborative work and also account for 
publication length.  While such distinctions might not be appropriate 
for other science disciplines, it is argued that they 
are instructive and necessary in mathematics. 
Such distinctions have been previously 
discussed within the mathematics community, though rarely (if ever) 
supported by the kinds of data provided here. 
Various statistical summaries which are 
presented in tables at the end of the paper are described 
in \S 4.  Some concluding remarks are given in \S 5. 

{\bf Note.} For the benefit of referees, more extensive data was submitted 
which will remain confidential. 

{\bf \S 2. Design of the study.} 
First a roster was created 
of mathematics faculty meeting the following criteria: 
($a$) doctorate earned in 2001 or earlier; ($b$) ranked faculty member 
in a mathematics department 
for the 2006-07 academic year at one of 
38 certain public universities (including the author's home 
institution Murray State University) 
with institutional profile, and presumably mission, 
comparable to Murray State's; ($c$) information concerning rank, doctorate school 
and year, as well as research interests and/or dissertation title 
readily available online; and ($d$) academic area in a mathematics field other 
than computer science, statistics, mathematics education, mathematics 
history, operations research, or actuarial science. 
This yielded a list of 366 individuals.  

These criteria are further explained/justified as follows.  
For ($a$), the data for a certain kind of research 
output collected for this study only covers publications 
up through calendar year 2006.  
Some of the summary statistics presented in this report are for 
certain five-year periods.  
Individuals with a doctorate earned in 2001 or earlier would 
therefore have had the opportunity 
to contribute for a minimum of five consecutive calendar years.  
For ($b$), only ranked faculty are considered since their typical 
professional responsibilities are 
not necessarily fully shared by other university faculty (lecturers, 
adjuncts, etc).  
Since this study of research productivity at regional 
public universities originated at Murray State, it made sense to look 
at other similar institutions.  
The 38 universities 
include all 15 of the public universities whose 
2005 Carnegie classification\myfootnote{See 
{\tt http://www.carnegiefoundation.org}} is Master's L and which are MSU 
benchmark schools (see \cite{MSUBenchmark})  
plus all 23 remaining public 
Master's L 
universities with mathematics departments in Kentucky 
and the Kentucky-area states of Arkansas, 
Indiana, Illinois, Missouri, 
North Carolina, Ohio, Tennessee, Virginia, and 
West Virginia. 
These 38 universities are nearly 
one-quarter of the 157 U.S.\ public Master's L universities with 
mathematics departments offering a traditional 
mathematics major.  (There are 166 U.S.\ public Master's L 
universities.)   
For ($c$), since one of the goals at the outset was to demonstrate the 
wide availability of data to support such a study, the focus was 
limited to information that could be obtained from  
highly accessible resources. 
Thus, roster information was 
obtained from open sources publicly available online: 
university and department webpages, online 
bulletins and library 
catalogs,\setcounter{myfn}{1}\renewcommand{\thefootnote}{\fnsymbol{myfn}} the 
Mathematics Genealogy 
Project\footnote{See 
{\tt http://www.genealogy.ams.org}}, 
and Google searches. 
For ($d$), the main reason for excluding individuals in these 
particular areas is that such faculty often have priorities for 
productivity other than research leading to the kinds of publications 
considered here. 
From the 38 departments 
there were 341 ranked mathematics faculty who were not eligible for 
this  roster: 15 computer 
scientists (9 with doctorate before 2002), 94 statisticians (73 with 
doctorate before 2002), 113 mathematics education specialists (88 
with doctorate before 2002), 3 mathematics historians (2 with 
doctorate before 2002), 12 operations researchers (10 with doctorate 
before 2002), 2 actuarial scientists (both with doctorate before 2002), 
74 others with doctorate between the years 
2002 and 2006 inclusive, 16 without a doctoral degree, and 12 others  
with insufficient information available online. 

\setcounter{myfn}{2} 
\renewcommand{\thefootnote}{\fnsymbol{myfn}}
Second, publication data for all faculty on the roster was obtained 
from 
MathSciNet.\footnote{See {\tt http://www.ams.org/mathscinet}\ \ (This 
subscription service is available 
on many university campuses.)}  
MathSciNet is a comprehensive 
searchable online database maintained by the American Mathematical 
Society (AMS) and which 
includes bibliographic information for mathematics 
papers published worldwide for the past 65+ years.  
MathSciNet has its origins as the monthly print publication 
Mathematical Reviews, begun in 1940 
with the goal of 
providing timely 
information on new contributions 
to mathematics research appearing in the literature. 
MathSciNet's reputation as one of the world's premiere databases of 
information 
on\setcounter{myfn}{3}\renewcommand{\thefootnote}{\fnsymbol{myfn}} mathematics 
literature owes to its reliability, wide availability, 
and vast coverage.\myfootnote{The only other comparable database is 
{\em Zentralblatt Math} at 
{\tt http://www.zentralblatt-math.org/zmath/en}} For example, 
each year all papers from nearly 750 mathematics journals 
worldwide have individual reviews in addition to bibliographic 
information added to the MathSciNet database.  
Bibliographic data is also indexed for 
hundreds of other journals related to the mathematical sciences 
which are not reviewed cover-to-cover. 
Many conference proceedings and other compendia of research papers are 
indexed as well. 
In what follows, the 
notation {\bf `MSN'} refers to papers 
with bibliographic information appearing on MathSciNet. 
It should be kept in mind that MSN publications represent only one 
kind of scholarly mathematical contribution. 
(For instance, some statistics and many education-related journals 
are not indexed by MathSciNet.) 
But in view of the position of the AMS  
that the legacy of the mathematics community is its publication 
record \cite{AMS}, MSN publications surely represent a very important kind of 
scholarly mathematical contribution.  
 
Third, a database 
was created by entering the raw bibliographic information on MSN 
publications for faculty on the roster into 
a text file.  This raw data was then compiled into statistical 
summaries using programs written in the 
computer algebra system Maple.   The biographical details obtained 
for the 366 individual faculty (complete name information, doctorate year, 
research area, university affiliation, etc) helped identify the MSN 
publications for each via MathSciNet's author search query. 
 
{\bf \S 3. Metrics considered.} 
This report will focus on statistics based on 
certain {\em paper} and {\em page} counts.  
These counts are tallied for each of the 366 individual faculty on 
the roster described in \S 2.  Certain averages for the faculty at each of 
the 38 universities are also calculated.  The counts are derived from 
MSN publication data.  
Since this study aims to focus on post-doctorate productivity, 
only MSN papers published \underline{after} the doctorate  
year are considered.  
Those items which are readily identified as errata, addenda, 
surveys/expositions, 
or research announcements (offering results without proofs)  
are excluded from paper-counting metrics. 
Pages from surveys/expositions  
and research announcements (offering results without proofs) 
are excluded from 
page-counting metrics. 
These 
exclusions amount to 
only a small fraction of the `attributable' papers 
(32.4 out of 2,156.3) 
and `attributable' pages 
(340.5 out of 26,386.6) 
analyzed in the study. 
Here and throughout, 
the adjective {\bf `attributable'} refers to 
the proportion of a paper obtained by dividing by the number of 
coauthors.  

In the category of 
papers, the metrics considered are single-authored papers, 
collaborative papers, attributable collaborative papers, and attributable 
papers. The latter is the sum of the single-authored and attributable 
collaborative papers.  For collaborative papers in mathematics, 
dividing by the number of coauthors is appropriate at times, for the 
following reasons. 
First, while it is not likely that contributions of all authors on a 
collaborative paper are exactly the same in terms of 
generating ideas, obtaining results, writing, etc, it is nonetheless 
reasonable to assume that each of the authors has made a significant 
contribution.  This is expressed by the AMS in its 
ethical guidelines for coauthorship, which state 
that all of the authors listed on a collaborative 
paper ``must have made a significant 
contribution to its content'' \cite{AMS}. 
Moreover, 
in mathematics the prevailing culture is to list authors 
alphabetically, so  
the bibliographic data usually makes no distinctions 
concerning coauthors' respective contributions.  
This practice reflects the facts that the typical end-product of a 
mathematical investigation is a new theorem or proof and that the 
relative merits of the input ideas contributed by collaborating 
researchers toward such an end-product can be difficult to 
distinguish. 
So, dividing by the number of coauthors is a workably equitable principle for 
accounting for an individual's relative contribution to a given collaborative 
paper.  
Second, 
from the point of view of editors or referees, there is no 
distinction made in standards for single-authored or coauthored 
papers.  
The evident principle is that the academic merit of a paper is 
independent of the number of authors. 
So in assessments of individual research productivity, 
the author of the single-authored paper can be unfairly 
disadvantaged when coauthors -- who can divide the labor of 
production amongst themselves -- receive the same individual credit as 
the single author.  
Further, equating the efforts of 
single authors and collaborative authors 
not only 
disadvantages the single author in this way but 
can be seen as effectively (and dubiously) giving credit to 
an author of a coauthored paper 
for the work of his or her collaborators. 
Thus, as a metric for assessing individual productivity in mathematics 
research, 
the total of items on the publication list is flawed in 
its conflation of contributions from single-authored and collaborative 
papers.  
(The notion that certain bibliometric indicators of scientific research 
performance  should account for coauthorship by 
dividing by the number of authors is not unusual, see  
for example in \cite{BatistaScientometrics}, \cite{Batista}, 
\cite{Gauffriau}, \cite{Larsen}.) 
On the other hand, collaborations can lead to research results that 
might not otherwise have been obtained and can say something very 
positive about a researcher or a research 
program, so it seems imperative that some metrics should 
specifically recognize this kind of contribution. 
Overall, as indicators of mathematics research productivity, both 
single-authored and collaborative publications are important. 

In the category of pages, the metrics considered are single-authored pages, 
collaborative pages, attributable collaborative pages, and attributable pages.  
The latter is the sum of the single-authored and attributable 
collaborative pages. 
In mathematics, such 
page counts are legitimate and in certain respects better 
measures than paper counts for assessing research 
output. 
Length and 
content for mathematics papers are reasonably viewed as 
proportional.  
In part this is because 
the recognized standards for journal writing encourage economy and 
demand that content represent new additions to the literature.  
Further, attempting to inflate page counts is risky since longer 
papers demand more of editors and particularly referees and 
take up more journal space, which is often at a premium.  
It can be argued that the paper is a far more arbitrary unit of volume 
than the page and far more susceptible to authors' stylistic choices.  
Moreover, a system of assessment in which 
the length of the publication list is the only metric 
that counts can create incentive for pursuing 
stratagems that have little or no academic 
merit, such as arbitrarily subdividing a paper to get multiple 
submissions.  
Some of that incentive is dispelled when 
page counts are also considered. 
One possible drawback of page-counting metrics is the implicit 
assumption that the content on any two given pages is roughly comparable, 
not only within a paper but also between 
papers in the same or even different journals. 
However, this flaw seems to be no worse than the assumption that the 
content of two {\em papers} is roughly comparable. 
Indeed, the CMS has 
recognized that as a ``pure, non-evaluative [i.e.\ objective] output measure'', 
a count of pages is not only legitimate but 
also less crude than a count of papers; further, counting pages as a 
measure of output volume can be 
seen as analogous to counting the 
monetary totals of grants and not just the number of grants won when 
assessing the input volume of research funding \cite{AppendixB}.  


{\bf \S 4. Comments on the statistical summaries.} 
In the first set of tables (Tables 
1.1--1.32), 
summary data is given for an octet of metrics 
(single-authored papers and pages, collaborative papers and pages, 
attributable collaborative papers and pages, and attributable 
papers and pages) for the 366-member faculty 
roster. Publication productivity is considered for the following time 
frames: The calendar years 2002-2006 (inclusive), 
the best five years, 
the best up-to-ten-year period, 
and total career productivity through calendar 
year 2006.  
For a given metric, 
the best five years for an individual are his/her best five consecutive 
calendar years during the career span from after the doctorate 
year up through 2006.  
The best up-to-ten-year period is  
the best ten consecutive years during the career span 
or the period from after  
the doctorate year through 2006, 
if the latter is less than ten years. 
In the language of \cite{Gauffriau} (see also \cite{Larsen}), the 
object of study in these tables 
is the 366-member faculty roster, the basic units 
are\setcounter{myfn}{1}\renewcommand{\thefootnote}{\fnsymbol{myfn}}
authors, and credits are attributed using complete counting 
(single-authored papers/pages, 
collaborative papers/pages) and 
complete-normalized\myfootnote{In the context of the metrics used 
here, ``normalizing'' refers to dividing by the number of coauthors.  
So, crediting using 
single-authored papers/pages could be viewed as both complete 
and complete-normalized.} counting (attributable collaborative 
papers/pages, attributable papers/pages).  

The second set of tables (Tables 2.1--2.16) gives data for 
faculty at each of the 38 universities 
which were part of this study.  
For each school and each of the eight metrics,  
the average per faculty member is given as well as the median 
and the average of the middle 50\%.  
These are presented for the four different time frames (2002-2006, 
best five years, 
best up-to-ten years, 
and career).  This data could be viewed as a companion to 
the publication data reported for faculty at research  
universities in Appendix L of \cite{NRC}.   
In the language of \cite{Gauffriau} (see also \cite{Larsen}), the 
objects of study in these tables 
are the 38 sets of faculty obtained by grouping the individuals from the 
366-member roster according to their 2006-07 university affiliation.  
{\em That is, the objects of study are not the universities 
themselves.} 
The basic units 
are authors, and credits are attributed using complete counting 
(single-authored papers/pages, 
collaborative papers/pages) 
and complete-normalized$^{*}$ 
counting (attributable collaborative 
papers/pages, attributable papers/pages).  
For the ``Totals'' data reported in these tables, the object of study 
is the entire 366-member faculty roster. 

Data in the first set of tables 
is reported in ``half-deciles'' on a 100-place scale (i.e.\ 100th 
place, 95th place, 90th place, etc).  
The 100th place is the highest mark in a list of 366.  The 0th place is the 
lowest mark out of 366. The 50th place is the median.  
The 65th place (for example) is computed as 
$.75 \times$ (238th highest mark) + $.25 \times$ (239th highest mark), 
corresponding to position $238.25 = 1 + 
65\cdot\left(\frac{366-1}{100}\right)$ in 
an ordered list of 366 marks. 
For the five-year and up-to-ten-year 
time frames in the first and second sets of data, 
biennial rather than annual averages are used 
for the reason that a biennium is 
more closely attuned to the natural life-cycle 
of a mathematics paper from initial idea to submission.
To compute a biennial average, divide an individual 
total by half the number of calendar years.   

The presentation of summary statistics begins on the next page. 

{\bf \S 5. Concluding remarks.} Distinguishing 
between single-authored and collaborative publications and accounting 
for publication length are necessary distinctions for informative 
statistical summaries of mathematics research output.  
It is difficult for single numerical measures (such as length of the 
publication list) to adequately draw such 
distinctions.  In addition to making such distinctions, 
the statistical summaries  
presented in this report use certain academic-biographical 
information to account for faculty research areas and to 
give a view of productivity for various time frames, 
most notably including the best five and up-to-ten years for each individual 
relative to each of the eight metrics used.  
It is hoped that this data will help support constructive discussions 
of the use of metrics in assessing mathematics research output, 
particularly at U.S.\ regional public universities. 
Further refinements of 
this data are planned for future reports, accounting for factors such as 
subject area within mathematics (using the AMS' 
Mathematics\setcounter{myfn}{2}\renewcommand{\thefootnote}{\fnsymbol{myfn}}
Subject 
Classification\myfootnote{See {\tt http://www.ams.org/msc}} 
scheme), diversity in the selection of 
journals in which an individual's work appears, diversity of 
coauthors, whether publications have appeared in archival journals or 
in proceedings or other research compendia, etc.  

{\bf Acknowledgment.} I thank my Murray State colleagues 
Wayne Bell and Andy Kellie 
for their encouragement and thoughtful feedback.

\newpage 
{\small 
\begin{center} 
\fbox{Tables 1.1 -- 1.8: time frame = 2002-2006} 

\ClarifyingRemark

\  

Table 1.1: Biennial averages for single-authored papers for the 
five-year period 2002-2006.

\begin{tabular}{|r|c||r|c||r|c||r|c||r|c||r|c||r|c|}
\hline 
100th & \parbox{0.31in}{\hfill \footnotesize 3.6} & 
85th & \parbox{0.31in}{\hfill \footnotesize 0.4} & 
70th & \parbox{0.31in}{\hfill \footnotesize 0} & 
55th & \parbox{0.31in}{\hfill \footnotesize 0} & 
40th & \parbox{0.31in}{\hfill \footnotesize 0} & 
25th & \parbox{0.31in}{\hfill \footnotesize 0} & 
10th & \parbox{0.31in}{\hfill \footnotesize 0}\\
\hline 
95th & \parbox{0.31in}{\hfill \footnotesize 1.2} & 
80th & \parbox{0.31in}{\hfill \footnotesize 0.4} & 
65th & \parbox{0.31in}{\hfill \footnotesize 0} & 
50th & \parbox{0.31in}{\hfill \footnotesize 0} & 
35th & \parbox{0.31in}{\hfill \footnotesize 0} & 
20th & \parbox{0.31in}{\hfill \footnotesize 0} & 
5th & \parbox{0.31in}{\hfill \footnotesize 0}\\
\hline 
90th & \parbox{0.31in}{\hfill \footnotesize 0.8} & 
75th & \parbox{0.31in}{\hfill \footnotesize 0} & 
60th & \parbox{0.31in}{\hfill \footnotesize 0} & 
45th & \parbox{0.31in}{\hfill \footnotesize 0} & 
30th & \parbox{0.31in}{\hfill \footnotesize 0} & 
15th & \parbox{0.31in}{\hfill \footnotesize 0} & 
0th & \parbox{0.31in}{\hfill \footnotesize 0}\\
\hline 
\end{tabular}

\ 

Table 1.2: Biennial averages for collaborative papers for the 
five-year period 2002-2006.

\begin{tabular}{|r|c||r|c||r|c||r|c||r|c||r|c||r|c|}
\hline 
100th & \parbox{0.31in}{\hfill \footnotesize 21.2} & 
85th & \parbox{0.31in}{\hfill \footnotesize 1.6} & 
70th & \parbox{0.31in}{\hfill \footnotesize 0.4} & 
55th & \parbox{0.31in}{\hfill \footnotesize 0} & 
40th & \parbox{0.31in}{\hfill \footnotesize 0} & 
25th & \parbox{0.31in}{\hfill \footnotesize 0} & 
10th & \parbox{0.31in}{\hfill \footnotesize 0}\\
\hline 
95th & \parbox{0.31in}{\hfill \footnotesize 2.8} & 
80th & \parbox{0.31in}{\hfill \footnotesize 1.2} & 
65th & \parbox{0.31in}{\hfill \footnotesize 0.4} & 
50th & \parbox{0.31in}{\hfill \footnotesize 0} & 
35th & \parbox{0.31in}{\hfill \footnotesize 0} & 
20th & \parbox{0.31in}{\hfill \footnotesize 0} & 
5th & \parbox{0.31in}{\hfill \footnotesize 0}\\
\hline 
90th & \parbox{0.31in}{\hfill \footnotesize 2.0} & 
75th & \parbox{0.31in}{\hfill \footnotesize 0.8} & 
60th & \parbox{0.31in}{\hfill \footnotesize 0.4} & 
45th & \parbox{0.31in}{\hfill \footnotesize 0} & 
30th & \parbox{0.31in}{\hfill \footnotesize 0} & 
15th & \parbox{0.31in}{\hfill \footnotesize 0} & 
0th & \parbox{0.31in}{\hfill \footnotesize 0}\\
\hline 
\end{tabular}

\ 

Table 1.3: Biennial averages for attributable collaborative papers for the 
five-year period 2002-2006.

\begin{tabular}{|r|c||r|c||r|c||r|c||r|c||r|c||r|c|}
\hline 
100th & \parbox{0.31in}{\hfill \footnotesize 8.53} & 
85th & \parbox{0.31in}{\hfill \footnotesize 0.60} & 
70th & \parbox{0.31in}{\hfill \footnotesize 0.20} & 
55th & \parbox{0.31in}{\hfill \footnotesize 0} & 
40th & \parbox{0.31in}{\hfill \footnotesize 0} & 
25th & \parbox{0.31in}{\hfill \footnotesize 0} & 
10th & \parbox{0.31in}{\hfill \footnotesize 0}\\
\hline 
95th & \parbox{0.31in}{\hfill \footnotesize 1.12} & 
80th & \parbox{0.31in}{\hfill \footnotesize 0.40} & 
65th & \parbox{0.31in}{\hfill \footnotesize 0.13} & 
50th & \parbox{0.31in}{\hfill \footnotesize 0} & 
35th & \parbox{0.31in}{\hfill \footnotesize 0} & 
20th & \parbox{0.31in}{\hfill \footnotesize 0} & 
5th & \parbox{0.31in}{\hfill \footnotesize 0}\\
\hline 
90th & \parbox{0.31in}{\hfill \footnotesize 0.80} & 
75th & \parbox{0.31in}{\hfill \footnotesize 0.33} & 
60th & \parbox{0.31in}{\hfill \footnotesize 0.13} & 
45th & \parbox{0.31in}{\hfill \footnotesize 0} & 
30th & \parbox{0.31in}{\hfill \footnotesize 0} & 
15th & \parbox{0.31in}{\hfill \footnotesize 0} & 
0th & \parbox{0.31in}{\hfill \footnotesize 0}\\
\hline 
\end{tabular}

\ 

Table 1.4: Biennial averages for attributable papers for the 
five-year period 2002-2006.

\begin{tabular}{|r|c||r|c||r|c||r|c||r|c||r|c||r|c|}
\hline 
100th & \parbox{0.31in}{\hfill \footnotesize 8.53} & 
85th & \parbox{0.31in}{\hfill \footnotesize 1.07} & 
70th & \parbox{0.31in}{\hfill \footnotesize 0.45} & 
55th & \parbox{0.31in}{\hfill \footnotesize 0.20} & 
40th & \parbox{0.31in}{\hfill \footnotesize 0} & 
25th & \parbox{0.31in}{\hfill \footnotesize 0} & 
10th & \parbox{0.31in}{\hfill \footnotesize 0}\\
\hline 
95th & \parbox{0.31in}{\hfill \footnotesize 1.90} & 
80th & \parbox{0.31in}{\hfill \footnotesize 0.80} & 
65th & \parbox{0.31in}{\hfill \footnotesize 0.40} & 
50th & \parbox{0.31in}{\hfill \footnotesize 0.13} & 
35th & \parbox{0.31in}{\hfill \footnotesize 0} & 
20th & \parbox{0.31in}{\hfill \footnotesize 0} & 
5th & \parbox{0.31in}{\hfill \footnotesize 0}\\
\hline 
90th & \parbox{0.31in}{\hfill \footnotesize 1.53} & 
75th & \parbox{0.31in}{\hfill \footnotesize 0.62} & 
60th & \parbox{0.31in}{\hfill \footnotesize 0.33} & 
45th & \parbox{0.31in}{\hfill \footnotesize 0} & 
30th & \parbox{0.31in}{\hfill \footnotesize 0} & 
15th & \parbox{0.31in}{\hfill \footnotesize 0} & 
0th & \parbox{0.31in}{\hfill \footnotesize 0}\\
\hline 
\end{tabular}

\ 

\underline{\hspace*{4in}} 

\

Table 1.5: Biennial averages for single-authored pages for the 
five-year period 2002-2006.

\begin{tabular}{|r|c||r|c||r|c||r|c||r|c||r|c||r|c|}
\hline 
100th & \parbox{0.31in}{\hfill \footnotesize 64.0} & 
85th & \parbox{0.31in}{\hfill \footnotesize 6.1} & 
70th & \parbox{0.31in}{\hfill \footnotesize 0} & 
55th & \parbox{0.31in}{\hfill \footnotesize 0} & 
40th & \parbox{0.31in}{\hfill \footnotesize 0} & 
25th & \parbox{0.31in}{\hfill \footnotesize 0} & 
10th & \parbox{0.31in}{\hfill \footnotesize 0}\\
\hline 
95th & \parbox{0.31in}{\hfill \footnotesize 16.2} & 
80th & \parbox{0.31in}{\hfill \footnotesize 3.6} & 
65th & \parbox{0.31in}{\hfill \footnotesize 0} & 
50th & \parbox{0.31in}{\hfill \footnotesize 0} & 
35th & \parbox{0.31in}{\hfill \footnotesize 0} & 
20th & \parbox{0.31in}{\hfill \footnotesize 0} & 
5th & \parbox{0.31in}{\hfill \footnotesize 0}\\
\hline 
90th & \parbox{0.31in}{\hfill \footnotesize 8.4} & 
75th & \parbox{0.31in}{\hfill \footnotesize 0} & 
60th & \parbox{0.31in}{\hfill \footnotesize 0} & 
45th & \parbox{0.31in}{\hfill \footnotesize 0} & 
30th & \parbox{0.31in}{\hfill \footnotesize 0} & 
15th & \parbox{0.31in}{\hfill \footnotesize 0} & 
0th & \parbox{0.31in}{\hfill \footnotesize 0}\\
\hline 
\end{tabular}

\ 

Table 1.6: Biennial averages for collaborative pages for the 
five-year period 2002-2006.

\begin{tabular}{|r|c||r|c||r|c||r|c||r|c||r|c||r|c|}
\hline 
100th & \parbox{0.31in}{\hfill \footnotesize 294.4} & 
85th & \parbox{0.31in}{\hfill \footnotesize 19.5} & 
70th & \parbox{0.31in}{\hfill \footnotesize 6.4} & 
55th & \parbox{0.31in}{\hfill \footnotesize 0} & 
40th & \parbox{0.31in}{\hfill \footnotesize 0} & 
25th & \parbox{0.31in}{\hfill \footnotesize 0} & 
10th & \parbox{0.31in}{\hfill \footnotesize 0}\\
\hline 
95th & \parbox{0.31in}{\hfill \footnotesize 46.2} & 
80th & \parbox{0.31in}{\hfill \footnotesize 13.2} & 
65th & \parbox{0.31in}{\hfill \footnotesize 4.8} & 
50th & \parbox{0.31in}{\hfill \footnotesize 0} & 
35th & \parbox{0.31in}{\hfill \footnotesize 0} & 
20th & \parbox{0.31in}{\hfill \footnotesize 0} & 
5th & \parbox{0.31in}{\hfill \footnotesize 0}\\
\hline 
90th & \parbox{0.31in}{\hfill \footnotesize 30.8} & 
75th & \parbox{0.31in}{\hfill \footnotesize 10.8} & 
60th & \parbox{0.31in}{\hfill \footnotesize 2.8} & 
45th & \parbox{0.31in}{\hfill \footnotesize 0} & 
30th & \parbox{0.31in}{\hfill \footnotesize 0} & 
15th & \parbox{0.31in}{\hfill \footnotesize 0} & 
0th & \parbox{0.31in}{\hfill \footnotesize 0}\\
\hline 
\end{tabular}

\ 

Table 1.7: Biennial averages for attributable collaborative pages for the 
five-year period 2002-2006.

\begin{tabular}{|r|c||r|c||r|c||r|c||r|c||r|c||r|c|}
\hline 
100th & {\footnotesize 119.03} & 
85th & \parbox{0.31in}{\hfill \footnotesize 8.63} & 
70th & \parbox{0.31in}{\hfill \footnotesize 2.89} & 
55th & \parbox{0.31in}{\hfill \footnotesize 0} & 
40th & \parbox{0.31in}{\hfill \footnotesize 0} & 
25th & \parbox{0.31in}{\hfill \footnotesize 0} & 
10th & \parbox{0.31in}{\hfill \footnotesize 0}\\
\hline 
95th & \parbox{0.31in}{\hfill \footnotesize 19.60} & 
80th & \parbox{0.31in}{\hfill \footnotesize 5.53} & 
65th & \parbox{0.31in}{\hfill \footnotesize 2.00} & 
50th & \parbox{0.31in}{\hfill \footnotesize 0} & 
35th & \parbox{0.31in}{\hfill \footnotesize 0} & 
20th & \parbox{0.31in}{\hfill \footnotesize 0} & 
5th & \parbox{0.31in}{\hfill \footnotesize 0}\\
\hline 
90th & \parbox{0.31in}{\hfill \footnotesize 12.25} & 
75th & \parbox{0.31in}{\hfill \footnotesize 4.15} & 
60th & \parbox{0.31in}{\hfill \footnotesize 0.93} & 
45th & \parbox{0.31in}{\hfill \footnotesize 0} & 
30th & \parbox{0.31in}{\hfill \footnotesize 0} & 
15th & \parbox{0.31in}{\hfill \footnotesize 0} & 
0th & \parbox{0.31in}{\hfill \footnotesize 0}\\
\hline 
\end{tabular}

\ 

Table 1.8: Biennial averages for attributable pages for the 
five-year period 2002-2006.

\begin{tabular}{|r|c||r|c||r|c||r|c||r|c||r|c||r|c|}
\hline 
100th & {\footnotesize 119.03} & 
85th & \parbox{0.31in}{\hfill \footnotesize 16.60} & 
70th & \parbox{0.31in}{\hfill \footnotesize 5.77} & 
55th & \parbox{0.31in}{\hfill \footnotesize 2.00} & 
40th & \parbox{0.31in}{\hfill \footnotesize 0} & 
25th & \parbox{0.31in}{\hfill \footnotesize 0} & 
10th & \parbox{0.31in}{\hfill \footnotesize 0}\\
\hline 
95th & \parbox{0.31in}{\hfill \footnotesize 30.26} & 
80th & \parbox{0.31in}{\hfill \footnotesize 11.73} & 
65th & \parbox{0.31in}{\hfill \footnotesize 4.50} & 
50th & \parbox{0.31in}{\hfill \footnotesize 0.93} & 
35th & \parbox{0.31in}{\hfill \footnotesize 0} & 
20th & \parbox{0.31in}{\hfill \footnotesize 0} & 
5th & \parbox{0.31in}{\hfill \footnotesize 0}\\
\hline 
90th & \parbox{0.31in}{\hfill \footnotesize 22.07} & 
75th & \parbox{0.31in}{\hfill \footnotesize 9.04} & 
60th & \parbox{0.31in}{\hfill \footnotesize 3.20} & 
45th & \parbox{0.31in}{\hfill \footnotesize 0} & 
30th & \parbox{0.31in}{\hfill \footnotesize 0} & 
15th & \parbox{0.31in}{\hfill \footnotesize 0} & 
0th & \parbox{0.31in}{\hfill \footnotesize 0}\\
\hline 
\end{tabular}
\end{center} 
}

\newpage 
{\small 
\begin{center} 
\fbox{Tables 1.9 -- 1.16: time frame = best five years}
 
\ClarifyingRemark
 
\ 

Table 1.9: Biennial averages for single-authored papers for the 
best five years. 

\begin{tabular}{|r|c||r|c||r|c||r|c||r|c||r|c||r|c|}
\hline 
100th & \parbox{0.31in}{\hfill \footnotesize 9.6} & 
85th & \parbox{0.31in}{\hfill \footnotesize 1.6} & 
70th & \parbox{0.31in}{\hfill \footnotesize 0.8} & 
55th & \parbox{0.31in}{\hfill \footnotesize 0.4} & 
40th & \parbox{0.31in}{\hfill \footnotesize 0.4} & 
25th & \parbox{0.31in}{\hfill \footnotesize 0} & 
10th & \parbox{0.31in}{\hfill \footnotesize 0}\\
\hline 
95th & \parbox{0.31in}{\hfill \footnotesize 2.7} & 
80th & \parbox{0.31in}{\hfill \footnotesize 1.2} & 
65th & \parbox{0.31in}{\hfill \footnotesize 0.8} & 
50th & \parbox{0.31in}{\hfill \footnotesize 0.4} & 
35th & \parbox{0.31in}{\hfill \footnotesize 0.4} & 
20th & \parbox{0.31in}{\hfill \footnotesize 0} & 
5th & \parbox{0.31in}{\hfill \footnotesize 0}\\
\hline 
90th & \parbox{0.31in}{\hfill \footnotesize 2.0} & 
75th & \parbox{0.31in}{\hfill \footnotesize 1.2} & 
60th & \parbox{0.31in}{\hfill \footnotesize 0.8} & 
45th & \parbox{0.31in}{\hfill \footnotesize 0.4} & 
30th & \parbox{0.31in}{\hfill \footnotesize 0} & 
15th & \parbox{0.31in}{\hfill \footnotesize 0} & 
0th & \parbox{0.31in}{\hfill \footnotesize 0}\\
\hline 
\end{tabular}

\ 

Table 1.10: Biennial averages for collaborative papers for the 
best five years. 

\begin{tabular}{|r|c||r|c||r|c||r|c||r|c||r|c||r|c|}
\hline 
100th & \parbox{0.31in}{\hfill \footnotesize 22.0} & 
85th & \parbox{0.31in}{\hfill \footnotesize 2.4} & 
70th & \parbox{0.31in}{\hfill \footnotesize 1.6} & 
55th & \parbox{0.31in}{\hfill \footnotesize 0.8} & 
40th & \parbox{0.31in}{\hfill \footnotesize 0.4} & 
25th & \parbox{0.31in}{\hfill \footnotesize 0} & 
10th & \parbox{0.31in}{\hfill \footnotesize 0}\\
\hline 
95th & \parbox{0.31in}{\hfill \footnotesize 4.8} & 
80th & \parbox{0.31in}{\hfill \footnotesize 2.0} & 
65th & \parbox{0.31in}{\hfill \footnotesize 1.2} & 
50th & \parbox{0.31in}{\hfill \footnotesize 0.4} & 
35th & \parbox{0.31in}{\hfill \footnotesize 0.4} & 
20th & \parbox{0.31in}{\hfill \footnotesize 0} & 
5th & \parbox{0.31in}{\hfill \footnotesize 0}\\
\hline 
90th & \parbox{0.31in}{\hfill \footnotesize 3.2} & 
75th & \parbox{0.31in}{\hfill \footnotesize 1.6} & 
60th & \parbox{0.31in}{\hfill \footnotesize 0.8} & 
45th & \parbox{0.31in}{\hfill \footnotesize 0.4} & 
30th & \parbox{0.31in}{\hfill \footnotesize 0.4} & 
15th & \parbox{0.31in}{\hfill \footnotesize 0} & 
0th & \parbox{0.31in}{\hfill \footnotesize 0}\\
\hline 
\end{tabular}

\ 

Table 1.11: Biennial averages for attributable collaborative papers for the 
best five years. 

\begin{tabular}{|r|c||r|c||r|c||r|c||r|c||r|c||r|c|}
\hline 
100th & \parbox{0.31in}{\hfill \footnotesize 8.70} & 
85th & \parbox{0.31in}{\hfill \footnotesize 1.13} & 
70th & \parbox{0.31in}{\hfill \footnotesize 0.65} & 
55th & \parbox{0.31in}{\hfill \footnotesize 0.33} & 
40th & \parbox{0.31in}{\hfill \footnotesize 0.20} & 
25th & \parbox{0.31in}{\hfill \footnotesize 0} & 
10th & \parbox{0.31in}{\hfill \footnotesize 0}\\
\hline 
95th & \parbox{0.31in}{\hfill \footnotesize 2.02} & 
80th & \parbox{0.31in}{\hfill \footnotesize 0.93} & 
65th & \parbox{0.31in}{\hfill \footnotesize 0.60} & 
50th & \parbox{0.31in}{\hfill \footnotesize 0.20} & 
35th & \parbox{0.31in}{\hfill \footnotesize 0.20} & 
20th & \parbox{0.31in}{\hfill \footnotesize 0} & 
5th & \parbox{0.31in}{\hfill \footnotesize 0}\\
\hline 
90th & \parbox{0.31in}{\hfill \footnotesize 1.40} & 
75th & \parbox{0.31in}{\hfill \footnotesize 0.79} & 
60th & \parbox{0.31in}{\hfill \footnotesize 0.40} & 
45th & \parbox{0.31in}{\hfill \footnotesize 0.20} & 
30th & \parbox{0.31in}{\hfill \footnotesize 0.13} & 
15th & \parbox{0.31in}{\hfill \footnotesize 0} & 
0th & \parbox{0.31in}{\hfill \footnotesize 0}\\
\hline 
\end{tabular}

\ 

Table 1.12: Biennial averages for attributable papers for the 
best five years. 

\begin{tabular}{|r|c||r|c||r|c||r|c||r|c||r|c||r|c|}
\hline 
100th & \parbox{0.31in}{\hfill \footnotesize 9.80} & 
85th & \parbox{0.31in}{\hfill \footnotesize 2.20} & 
70th & \parbox{0.31in}{\hfill \footnotesize 1.53} & 
55th & \parbox{0.31in}{\hfill \footnotesize 1.00} & 
40th & \parbox{0.31in}{\hfill \footnotesize 0.60} & 
25th & \parbox{0.31in}{\hfill \footnotesize 0.40} & 
10th & \parbox{0.31in}{\hfill \footnotesize 0}\\
\hline 
95th & \parbox{0.31in}{\hfill \footnotesize 3.52} & 
80th & \parbox{0.31in}{\hfill \footnotesize 1.93} & 
65th & \parbox{0.31in}{\hfill \footnotesize 1.22} & 
50th & \parbox{0.31in}{\hfill \footnotesize 0.80} & 
35th & \parbox{0.31in}{\hfill \footnotesize 0.47} & 
20th & \parbox{0.31in}{\hfill \footnotesize 0.20} & 
5th & \parbox{0.31in}{\hfill \footnotesize 0}\\
\hline 
90th & \parbox{0.31in}{\hfill \footnotesize 2.60} & 
75th & \parbox{0.31in}{\hfill \footnotesize 1.61} & 
60th & \parbox{0.31in}{\hfill \footnotesize 1.20} & 
45th & \parbox{0.31in}{\hfill \footnotesize 0.80} & 
30th & \parbox{0.31in}{\hfill \footnotesize 0.40} & 
15th & \parbox{0.31in}{\hfill \footnotesize 0.13} & 
0th & \parbox{0.31in}{\hfill \footnotesize 0}\\
\hline 
\end{tabular}

\ 

\underline{\hspace*{4in}} 

\

Table 1.13: Biennial averages for single-authored pages for the 
best five years. 

\begin{tabular}{|r|c||r|c||r|c||r|c||r|c||r|c||r|c|}
\hline 
100th & \parbox{0.31in}{\hfill \footnotesize 109.6} & 
85th & \parbox{0.31in}{\hfill \footnotesize 18.8} & 
70th & \parbox{0.31in}{\hfill \footnotesize 11.6} & 
55th & \parbox{0.31in}{\hfill \footnotesize 7.2} & 
40th & \parbox{0.31in}{\hfill \footnotesize 4.0} & 
25th & \parbox{0.31in}{\hfill \footnotesize 0} & 
10th & \parbox{0.31in}{\hfill \footnotesize 0}\\
\hline 
95th & \parbox{0.31in}{\hfill \footnotesize 32.0} & 
80th & \parbox{0.31in}{\hfill \footnotesize 16.0} & 
65th & \parbox{0.31in}{\hfill \footnotesize 9.6} & 
50th & \parbox{0.31in}{\hfill \footnotesize 6.0} & 
35th & \parbox{0.31in}{\hfill \footnotesize 2.0} & 
20th & \parbox{0.31in}{\hfill \footnotesize 0} & 
5th & \parbox{0.31in}{\hfill \footnotesize 0}\\
\hline 
90th & \parbox{0.31in}{\hfill \footnotesize 24.8} & 
75th & \parbox{0.31in}{\hfill \footnotesize 13.5} & 
60th & \parbox{0.31in}{\hfill \footnotesize 8.8} & 
45th & \parbox{0.31in}{\hfill \footnotesize 4.8} & 
30th & \parbox{0.31in}{\hfill \footnotesize 0} & 
15th & \parbox{0.31in}{\hfill \footnotesize 0} & 
0th & \parbox{0.31in}{\hfill \footnotesize 0}\\
\hline 
\end{tabular}

\ 

Table 1.14: Biennial averages for collaborative pages for the 
best five years. 

\begin{tabular}{|r|c||r|c||r|c||r|c||r|c||r|c||r|c|}
\hline 
100th & \parbox{0.31in}{\hfill \footnotesize 294.4} & 
85th & \parbox{0.31in}{\hfill \footnotesize 37.7} & 
70th & \parbox{0.31in}{\hfill \footnotesize 18.8} & 
55th & \parbox{0.31in}{\hfill \footnotesize 10.4} & 
40th & \parbox{0.31in}{\hfill \footnotesize 4.8} & 
25th & \parbox{0.31in}{\hfill \footnotesize 0} & 
10th & \parbox{0.31in}{\hfill \footnotesize 0}\\
\hline 
95th & \parbox{0.31in}{\hfill \footnotesize 69.8} & 
80th & \parbox{0.31in}{\hfill \footnotesize 29.2} & 
65th & \parbox{0.31in}{\hfill \footnotesize 15.3} & 
50th & \parbox{0.31in}{\hfill \footnotesize 8.4} & 
35th & \parbox{0.31in}{\hfill \footnotesize 3.2} & 
20th & \parbox{0.31in}{\hfill \footnotesize 0} & 
5th & \parbox{0.31in}{\hfill \footnotesize 0}\\
\hline 
90th & \parbox{0.31in}{\hfill \footnotesize 45.8} & 
75th & \parbox{0.31in}{\hfill \footnotesize 22.8} & 
60th & \parbox{0.31in}{\hfill \footnotesize 12.8} & 
45th & \parbox{0.31in}{\hfill \footnotesize 6.4} & 
30th & \parbox{0.31in}{\hfill \footnotesize 1.2} & 
15th & \parbox{0.31in}{\hfill \footnotesize 0} & 
0th & \parbox{0.31in}{\hfill \footnotesize 0}\\
\hline 
\end{tabular}

\ 

Table 1.15: Biennial averages for attributable collaborative pages for the 
best five years. 

\begin{tabular}{|r|c||r|c||r|c||r|c||r|c||r|c||r|c|}
\hline 
100th & {\footnotesize 119.03} & 
85th & \parbox{0.31in}{\hfill \footnotesize 14.95} & 
70th & \parbox{0.31in}{\hfill \footnotesize 7.77} & 
55th & \parbox{0.31in}{\hfill \footnotesize 4.40} & 
40th & \parbox{0.31in}{\hfill \footnotesize 2.13} & 
25th & \parbox{0.31in}{\hfill \footnotesize 0} & 
10th & \parbox{0.31in}{\hfill \footnotesize 0}\\
\hline 
95th & \hspace*{0.01in} \parbox{0.31in}{\hfill \footnotesize 27.24} & 
80th & \parbox{0.31in}{\hfill \footnotesize 12.87} & 
65th & \parbox{0.31in}{\hfill \footnotesize 6.65} & 
50th & \parbox{0.31in}{\hfill \footnotesize 3.57} & 
35th & \parbox{0.31in}{\hfill \footnotesize 1.40} & 
20th & \parbox{0.31in}{\hfill \footnotesize 0} & 
5th & \parbox{0.31in}{\hfill \footnotesize 0}\\
\hline 
90th & \hspace*{0.01in} \parbox{0.31in}{\hfill \footnotesize 20.62} & 
75th & \parbox{0.31in}{\hfill \footnotesize 10.15} & 
60th & \parbox{0.31in}{\hfill \footnotesize 5.33} & 
45th & \parbox{0.31in}{\hfill \footnotesize 2.80} & 
30th & \parbox{0.31in}{\hfill \footnotesize 0.60} & 
15th & \parbox{0.31in}{\hfill \footnotesize 0} & 
0th & \parbox{0.31in}{\hfill \footnotesize 0}\\
\hline 
\end{tabular}

\ 

Table 1.16: Biennial averages for attributable pages for the 
best five years. 

\begin{tabular}{|r|c||r|c||r|c||r|c||r|c||r|c||r|c|}
\hline 
100th & {\footnotesize 121.20} & 
85th & \parbox{0.31in}{\hfill \footnotesize 27.49} & 
70th & \parbox{0.31in}{\hfill \footnotesize 18.44} & 
55th & \parbox{0.31in}{\hfill \footnotesize 12.18} & 
40th & \parbox{0.31in}{\hfill \footnotesize 7.20} & 
25th & \parbox{0.31in}{\hfill \footnotesize 3.95} & 
10th & \parbox{0.31in}{\hfill \footnotesize 0}\\
\hline 
95th & \hspace*{0.01in} \parbox{0.31in}{\hfill \footnotesize 46.71} & 
80th & \parbox{0.31in}{\hfill \footnotesize 23.80} & 
65th & \parbox{0.31in}{\hfill \footnotesize 16.40} & 
50th & \parbox{0.31in}{\hfill \footnotesize 10.40} & 
35th & \parbox{0.31in}{\hfill \footnotesize 6.00} & 
20th & \parbox{0.31in}{\hfill \footnotesize 2.40} & 
5th & \parbox{0.31in}{\hfill \footnotesize 0}\\
\hline 
90th & \hspace*{0.01in} \parbox{0.31in}{\hfill \footnotesize 33.74} & 
75th & \parbox{0.31in}{\hfill \footnotesize 20.83} & 
60th & \parbox{0.31in}{\hfill \footnotesize 14.20} & 
45th & \parbox{0.31in}{\hfill \footnotesize 8.55} & 
30th & \parbox{0.31in}{\hfill \footnotesize 5.10} & 
15th & \parbox{0.31in}{\hfill \footnotesize 1.38} & 
0th & \parbox{0.31in}{\hfill \footnotesize 0}\\
\hline 
\end{tabular}
\end{center} 
}

\newpage 
{\small 
\begin{center} 
\fbox{Tables 1.17 -- 1.24: time frame = best up-to-ten years}
 
\ClarifyingRemark
 
\ 

Table 1.17: Biennial averages for single-authored papers for the 
best up-to-ten years. 

\begin{tabular}{|r|c||r|c||r|c||r|c||r|c||r|c||r|c|}
\hline 
100th & \parbox{0.31in}{\hfill \footnotesize 8.40} & 
85th & \parbox{0.31in}{\hfill \footnotesize 1.03} & 
70th & \parbox{0.31in}{\hfill \footnotesize 0.60} & 
55th & \parbox{0.31in}{\hfill \footnotesize 0.40} & 
40th & \parbox{0.31in}{\hfill \footnotesize 0.20} & 
25th & \parbox{0.31in}{\hfill \footnotesize 0} & 
10th & \parbox{0.31in}{\hfill \footnotesize 0}\\
\hline 
95th & \parbox{0.31in}{\hfill \footnotesize 1.80} & 
80th & \parbox{0.31in}{\hfill \footnotesize 0.80} & 
65th & \parbox{0.31in}{\hfill \footnotesize 0.58} & 
50th & \parbox{0.31in}{\hfill \footnotesize 0.31} & 
35th & \parbox{0.31in}{\hfill \footnotesize 0.20} & 
20th & \parbox{0.31in}{\hfill \footnotesize 0} & 
5th & \parbox{0.31in}{\hfill \footnotesize 0}\\
\hline 
90th & \parbox{0.31in}{\hfill \footnotesize 1.40} & 
75th & \parbox{0.31in}{\hfill \footnotesize 0.80} & 
60th & \parbox{0.31in}{\hfill \footnotesize 0.40} & 
45th & \parbox{0.31in}{\hfill \footnotesize 0.20} & 
30th & \parbox{0.31in}{\hfill \footnotesize 0} & 
15th & \parbox{0.31in}{\hfill \footnotesize 0} & 
0th & \parbox{0.31in}{\hfill \footnotesize 0}\\
\hline 
\end{tabular}

\ 

Table 1.18: Biennial averages for collaborative papers for the 
best up-to-ten years. 

\begin{tabular}{|r|c||r|c||r|c||r|c||r|c||r|c||r|c|}
\hline 
100th & \parbox{0.31in}{\hfill \footnotesize 21.40} & 
85th & \parbox{0.31in}{\hfill \footnotesize 2.00} & 
70th & \parbox{0.31in}{\hfill \footnotesize 1.00} & 
55th & \parbox{0.31in}{\hfill \footnotesize 0.40} & 
40th & \parbox{0.31in}{\hfill \footnotesize 0.20} & 
25th & \parbox{0.31in}{\hfill \footnotesize 0} & 
10th & \parbox{0.31in}{\hfill \footnotesize 0}\\
\hline 
95th & \parbox{0.31in}{\hfill \footnotesize 3.40} & 
80th & \parbox{0.31in}{\hfill \footnotesize 1.50} & 
65th & \parbox{0.31in}{\hfill \footnotesize 0.80} & 
50th & \parbox{0.31in}{\hfill \footnotesize 0.40} & 
35th & \parbox{0.31in}{\hfill \footnotesize 0.20} & 
20th & \parbox{0.31in}{\hfill \footnotesize 0} & 
5th & \parbox{0.31in}{\hfill \footnotesize 0}\\
\hline 
90th & \parbox{0.31in}{\hfill \footnotesize 2.40} & 
75th & \parbox{0.31in}{\hfill \footnotesize 1.20} & 
60th & \parbox{0.31in}{\hfill \footnotesize 0.60} & 
45th & \parbox{0.31in}{\hfill \footnotesize 0.21} & 
30th & \parbox{0.31in}{\hfill \footnotesize 0.20} & 
15th & \parbox{0.31in}{\hfill \footnotesize 0} & 
0th & \parbox{0.31in}{\hfill \footnotesize 0}\\
\hline 
\end{tabular}

\ 

Table 1.19: Biennial averages for attributable collaborative papers for the 
best up-to-ten years. 

\begin{tabular}{|r|c||r|c||r|c||r|c||r|c||r|c||r|c|}
\hline 
100th & \parbox{0.31in}{\hfill \footnotesize 8.53} & 
85th & \parbox{0.31in}{\hfill \footnotesize 0.80} & 
70th & \parbox{0.31in}{\hfill \footnotesize 0.43} & 
55th & \parbox{0.31in}{\hfill \footnotesize 0.20} & 
40th & \parbox{0.31in}{\hfill \footnotesize 0.10} & 
25th & \parbox{0.31in}{\hfill \footnotesize 0} & 
10th & \parbox{0.31in}{\hfill \footnotesize 0}\\
\hline 
95th & \parbox{0.31in}{\hfill \footnotesize 1.42} & 
80th & \parbox{0.31in}{\hfill \footnotesize 0.67} & 
65th & \parbox{0.31in}{\hfill \footnotesize 0.37} & 
50th & \parbox{0.31in}{\hfill \footnotesize 0.17} & 
35th & \parbox{0.31in}{\hfill \footnotesize 0.10} & 
20th & \parbox{0.31in}{\hfill \footnotesize 0} & 
5th & \parbox{0.31in}{\hfill \footnotesize 0}\\
\hline 
90th & \parbox{0.31in}{\hfill \footnotesize 1.08} & 
75th & \parbox{0.31in}{\hfill \footnotesize 0.50} & 
60th & \parbox{0.31in}{\hfill \footnotesize 0.30} & 
45th & \parbox{0.31in}{\hfill \footnotesize 0.10} & 
30th & \parbox{0.31in}{\hfill \footnotesize 0.07} & 
15th & \parbox{0.31in}{\hfill \footnotesize 0} & 
0th & \parbox{0.31in}{\hfill \footnotesize 0}\\
\hline 
\end{tabular}

\ 

Table 1.20: Biennial averages for attributable papers for the 
best up-to-ten years. 

\begin{tabular}{|r|c||r|c||r|c||r|c||r|c||r|c||r|c|}
\hline 
100th & \parbox{0.31in}{\hfill \footnotesize 9.10} & 
85th & \parbox{0.31in}{\hfill \footnotesize 1.74} & 
70th & \parbox{0.31in}{\hfill \footnotesize 1.08} & 
55th & \parbox{0.31in}{\hfill \footnotesize 0.67} & 
40th & \parbox{0.31in}{\hfill \footnotesize 0.40} & 
25th & \parbox{0.31in}{\hfill \footnotesize 0.20} & 
10th & \parbox{0.31in}{\hfill \footnotesize 0}\\
\hline 
95th & \parbox{0.31in}{\hfill \footnotesize 2.62} & 
80th & \parbox{0.31in}{\hfill \footnotesize 1.45} & 
65th & \parbox{0.31in}{\hfill \footnotesize 0.90} & 
50th & \parbox{0.31in}{\hfill \footnotesize 0.54} & 
35th & \parbox{0.31in}{\hfill \footnotesize 0.30} & 
20th & \parbox{0.31in}{\hfill \footnotesize 0.10} & 
5th & \parbox{0.31in}{\hfill \footnotesize 0}\\
\hline 
90th & \parbox{0.31in}{\hfill \footnotesize 2.09} & 
75th & \parbox{0.31in}{\hfill \footnotesize 1.20} & 
60th & \parbox{0.31in}{\hfill \footnotesize 0.77} & 
45th & \parbox{0.31in}{\hfill \footnotesize 0.49} & 
30th & \parbox{0.31in}{\hfill \footnotesize 0.20} & 
15th & \parbox{0.31in}{\hfill \footnotesize 0.07} & 
0th & \parbox{0.31in}{\hfill \footnotesize 0}\\
\hline 
\end{tabular}

\ 

\underline{\hspace*{4in}} 

\

Table 1.21: Biennial averages for single-authored pages for the 
best up-to-ten years. 

\begin{tabular}{|r|c||r|c||r|c||r|c||r|c||r|c||r|c|}
\hline 
100th & \parbox{0.31in}{\hfill \footnotesize 92.20} & 
85th & \parbox{0.31in}{\hfill \footnotesize 13.00} & 
70th & \parbox{0.31in}{\hfill \footnotesize 7.20} & 
55th & \parbox{0.31in}{\hfill \footnotesize 4.15} & 
40th & \parbox{0.31in}{\hfill \footnotesize 2.20} & 
25th & \parbox{0.31in}{\hfill \footnotesize 0} & 
10th & \parbox{0.31in}{\hfill \footnotesize 0}\\
\hline 
95th & \parbox{0.31in}{\hfill \footnotesize 23.50} & 
80th & \parbox{0.31in}{\hfill \footnotesize 10.40} & 
65th & \parbox{0.31in}{\hfill \footnotesize 5.85} & 
50th & \parbox{0.31in}{\hfill \footnotesize 3.20} & 
35th & \parbox{0.31in}{\hfill \footnotesize 1.00} & 
20th & \parbox{0.31in}{\hfill \footnotesize 0} & 
5th & \parbox{0.31in}{\hfill \footnotesize 0}\\
\hline 
90th & \parbox{0.31in}{\hfill \footnotesize 15.80} & 
75th & \parbox{0.31in}{\hfill \footnotesize 8.58} & 
60th & \parbox{0.31in}{\hfill \footnotesize 5.00} & 
45th & \parbox{0.31in}{\hfill \footnotesize 2.60} & 
30th & \parbox{0.31in}{\hfill \footnotesize 0} & 
15th & \parbox{0.31in}{\hfill \footnotesize 0} & 
0th & \parbox{0.31in}{\hfill \footnotesize 0}\\
\hline 
\end{tabular}

\ 

Table 1.22: Biennial averages for collaborative pages for the 
best up-to-ten years. 

\begin{tabular}{|r|c||r|c||r|c||r|c||r|c||r|c||r|c|}
\hline 
100th & {\footnotesize 277.80} & 
85th & \parbox{0.31in}{\hfill \footnotesize 24.95} & 
70th & \parbox{0.31in}{\hfill \footnotesize 12.70} & 
55th & \parbox{0.31in}{\hfill \footnotesize 6.24} & 
40th & \parbox{0.31in}{\hfill \footnotesize 2.60} & 
25th & \parbox{0.31in}{\hfill \footnotesize 0} & 
10th & \parbox{0.31in}{\hfill \footnotesize 0}\\
\hline 
95th & \hspace*{0.005in} \parbox{0.31in}{\hfill \footnotesize 48.70} & 
80th & \parbox{0.31in}{\hfill \footnotesize 20.89} & 
65th & \parbox{0.31in}{\hfill \footnotesize 10.00} & 
50th & \parbox{0.31in}{\hfill \footnotesize 4.78} & 
35th & \parbox{0.31in}{\hfill \footnotesize 1.75} & 
20th & \parbox{0.31in}{\hfill \footnotesize 0} & 
5th & \parbox{0.31in}{\hfill \footnotesize 0}\\
\hline 
90th & \hspace*{0.005in} \parbox{0.31in}{\hfill \footnotesize 33.60} & 
75th & \parbox{0.31in}{\hfill \footnotesize 15.35} & 
60th & \parbox{0.31in}{\hfill \footnotesize 8.40} & 
45th & \parbox{0.31in}{\hfill \footnotesize 3.57} & 
30th & \parbox{0.31in}{\hfill \footnotesize 0.60} & 
15th & \parbox{0.31in}{\hfill \footnotesize 0} & 
0th & \parbox{0.31in}{\hfill \footnotesize 0}\\
\hline 
\end{tabular}

\ 

Table 1.23: Biennial averages for attributable collaborative pages for the 
best up-to-ten years. 

\begin{tabular}{|r|c||r|c||r|c||r|c||r|c||r|c||r|c|}
\hline 
100th & {\footnotesize 111.77} & 
85th & \parbox{0.31in}{\hfill \footnotesize 11.23} & 
70th & \parbox{0.31in}{\hfill \footnotesize 5.07} & 
55th & \parbox{0.31in}{\hfill \footnotesize 2.58} & 
40th & \parbox{0.31in}{\hfill \footnotesize 1.20} & 
25th & \parbox{0.31in}{\hfill \footnotesize 0} & 
10th & \parbox{0.31in}{\hfill \footnotesize 0}\\
\hline 
95th & \hspace*{0.005in} \parbox{0.31in}{\hfill \footnotesize 19.79} & 
80th & \parbox{0.31in}{\hfill \footnotesize 9.13} & 
65th & \parbox{0.31in}{\hfill \footnotesize 4.32} & 
50th & \parbox{0.31in}{\hfill \footnotesize 2.00} & 
35th & \parbox{0.31in}{\hfill \footnotesize 0.80} & 
20th & \parbox{0.31in}{\hfill \footnotesize 0} & 
5th & \parbox{0.31in}{\hfill \footnotesize 0}\\
\hline 
90th & \hspace*{0.005in} \parbox{0.31in}{\hfill \footnotesize 14.64} & 
75th & \parbox{0.31in}{\hfill \footnotesize 6.86} & 
60th & \parbox{0.31in}{\hfill \footnotesize 3.48} & 
45th & \parbox{0.31in}{\hfill \footnotesize 1.67} & 
30th & \parbox{0.31in}{\hfill \footnotesize 0.30} & 
15th & \parbox{0.31in}{\hfill \footnotesize 0} & 
0th & \parbox{0.31in}{\hfill \footnotesize 0}\\
\hline 
\end{tabular}

\ 

Table 1.24: Biennial averages for attributable pages for the 
best up-to-ten years. 

\begin{tabular}{|r|c||r|c||r|c||r|c||r|c||r|c||r|c|}
\hline 
100th & {\footnotesize 111.77} & 
85th & \parbox{0.31in}{\hfill \footnotesize 21.25} & 
70th & \parbox{0.31in}{\hfill \footnotesize 12.33} & 
55th & \parbox{0.31in}{\hfill \footnotesize 7.79} & 
40th & \parbox{0.31in}{\hfill \footnotesize 4.60} & 
25th & \parbox{0.31in}{\hfill \footnotesize 2.11} & 
10th & \parbox{0.31in}{\hfill \footnotesize 0}\\
\hline 
95th & \hspace*{0.005in} \parbox{0.31in}{\hfill \footnotesize 33.95} & 
80th & \parbox{0.31in}{\hfill \footnotesize 17.57} & 
65th & \parbox{0.31in}{\hfill \footnotesize 10.97} & 
50th & \parbox{0.31in}{\hfill \footnotesize 6.60} & 
35th & \parbox{0.31in}{\hfill \footnotesize 3.71} & 
20th & \parbox{0.31in}{\hfill \footnotesize 1.30} & 
5th & \parbox{0.31in}{\hfill \footnotesize 0}\\
\hline 
90th & \hspace*{0.005in} \parbox{0.31in}{\hfill \footnotesize 26.21} & 
75th & \parbox{0.31in}{\hfill \footnotesize 15.19} & 
60th & \parbox{0.31in}{\hfill \footnotesize 9.20} & 
45th & \parbox{0.31in}{\hfill \footnotesize 5.74} & 
30th & \parbox{0.31in}{\hfill \footnotesize 2.80} & 
15th & \parbox{0.31in}{\hfill \footnotesize 0.69} & 
0th & \parbox{0.31in}{\hfill \footnotesize 0}\\
\hline 
\end{tabular}
\end{center} 
}

\newpage 
{\small 
\begin{center} 
\fbox{Tables 1.25 -- 1.32: time frame = career}

\ClarifyingRemark
 
\ 

Table 1.25: Career totals for single-authored papers. 

\begin{tabular}{|r|c||r|c||r|c||r|c||r|c||r|c||r|c|}
\hline 
100th & \parbox{0.31in}{\hfill \footnotesize 73} & 
85th & \parbox{0.31in}{\hfill \footnotesize 6} & 
70th & \parbox{0.31in}{\hfill \footnotesize 3} & 
55th & \parbox{0.31in}{\hfill \footnotesize 2} & 
40th & \parbox{0.31in}{\hfill \footnotesize 1} & 
25th & \parbox{0.31in}{\hfill \footnotesize 0} & 
10th & \parbox{0.31in}{\hfill \footnotesize 0}\\
\hline 
95th & \parbox{0.31in}{\hfill \footnotesize 12} & 
80th & \parbox{0.31in}{\hfill \footnotesize 4} & 
65th & \parbox{0.31in}{\hfill \footnotesize 3} & 
50th & \parbox{0.31in}{\hfill \footnotesize 1.5} & 
35th & \parbox{0.31in}{\hfill \footnotesize 1} & 
20th & \parbox{0.31in}{\hfill \footnotesize 0} & 
5th & \parbox{0.31in}{\hfill \footnotesize 0}\\
\hline 
90th & \parbox{0.31in}{\hfill \footnotesize 7.5} & 
75th & \parbox{0.31in}{\hfill \footnotesize 4} & 
60th & \parbox{0.31in}{\hfill \footnotesize 2} & 
45th & \parbox{0.31in}{\hfill \footnotesize 1} & 
30th & \parbox{0.31in}{\hfill \footnotesize 0} & 
15th & \parbox{0.31in}{\hfill \footnotesize 0} & 
0th & \parbox{0.31in}{\hfill \footnotesize 0}\\
\hline 
\end{tabular}

\ 

Table 1.26: Career totals for collaborative papers. 

\begin{tabular}{|r|c||r|c||r|c||r|c||r|c||r|c||r|c|}
\hline 
100th & \parbox{0.31in}{\hfill \footnotesize 188} & 
85th & \parbox{0.31in}{\hfill \footnotesize 11} & 
70th & \parbox{0.31in}{\hfill \footnotesize 5} & 
55th & \parbox{0.31in}{\hfill \footnotesize 2} & 
40th & \parbox{0.31in}{\hfill \footnotesize 1} & 
25th & \parbox{0.31in}{\hfill \footnotesize 0} & 
10th & \parbox{0.31in}{\hfill \footnotesize 0}\\
\hline 
95th & \parbox{0.31in}{\hfill \footnotesize 22} & 
80th & \parbox{0.31in}{\hfill \footnotesize 8} & 
65th & \parbox{0.31in}{\hfill \footnotesize 4} & 
50th & \parbox{0.31in}{\hfill \footnotesize 2} & 
35th & \parbox{0.31in}{\hfill \footnotesize 1} & 
20th & \parbox{0.31in}{\hfill \footnotesize 0} & 
5th & \parbox{0.31in}{\hfill \footnotesize 0}\\
\hline 
90th & \parbox{0.31in}{\hfill \footnotesize 16} & 
75th & \parbox{0.31in}{\hfill \footnotesize 7} & 
60th & \parbox{0.31in}{\hfill \footnotesize 3} & 
45th & \parbox{0.31in}{\hfill \footnotesize 1} & 
30th & \parbox{0.31in}{\hfill \footnotesize 1} & 
15th & \parbox{0.31in}{\hfill \footnotesize 0} & 
0th & \parbox{0.31in}{\hfill \footnotesize 0}\\
\hline 
\end{tabular}

\ 

Table 1.27: Career totals for attributable collaborative papers. 

\begin{tabular}{|r|c||r|c||r|c||r|c||r|c||r|c||r|c|}
\hline 
100th & \parbox{0.31in}{\hfill \footnotesize 74.75} & 
85th & \parbox{0.31in}{\hfill \footnotesize 4.83} & 
70th & \parbox{0.31in}{\hfill \footnotesize 2.50} & 
55th & \parbox{0.31in}{\hfill \footnotesize 1.00} & 
40th & \parbox{0.31in}{\hfill \footnotesize 0.50} & 
25th & \parbox{0.31in}{\hfill \footnotesize 0} & 
10th & \parbox{0.31in}{\hfill \footnotesize 0}\\
\hline 
95th & \parbox{0.31in}{\hfill \footnotesize 9.46} & 
80th & \parbox{0.31in}{\hfill \footnotesize 3.50} & 
65th & \parbox{0.31in}{\hfill \footnotesize 1.83} & 
50th & \parbox{0.31in}{\hfill \footnotesize 0.83} & 
35th & \parbox{0.31in}{\hfill \footnotesize 0.50} & 
20th & \parbox{0.31in}{\hfill \footnotesize 0} & 
5th & \parbox{0.31in}{\hfill \footnotesize 0}\\
\hline 
90th & \parbox{0.31in}{\hfill \footnotesize 6.50} & 
75th & \parbox{0.31in}{\hfill \footnotesize 2.83} & 
60th & \parbox{0.31in}{\hfill \footnotesize 1.50} & 
45th & \parbox{0.31in}{\hfill \footnotesize 0.50} & 
30th & \parbox{0.31in}{\hfill \footnotesize 0.33} & 
15th & \parbox{0.31in}{\hfill \footnotesize 0} & 
0th & \parbox{0.31in}{\hfill \footnotesize 0}\\
\hline 
\end{tabular}

\ 

Table 1.28: Career totals  for attributable papers. 

\begin{tabular}{|r|c||r|c||r|c||r|c||r|c||r|c||r|c|}
\hline 
100th & \parbox{0.31in}{\hfill \footnotesize 94.67} & 
85th & \parbox{0.31in}{\hfill \footnotesize 11.04} & 
70th & \parbox{0.31in}{\hfill \footnotesize 5.50} & 
55th & \parbox{0.31in}{\hfill \footnotesize 3.50} & 
40th & \parbox{0.31in}{\hfill \footnotesize 2.00} & 
25th & \parbox{0.31in}{\hfill \footnotesize 1.00} & 
10th & \parbox{0.31in}{\hfill \footnotesize 0}\\
\hline 
95th & \parbox{0.31in}{\hfill \footnotesize 20.42} & 
80th & \parbox{0.31in}{\hfill \footnotesize 8.25} & 
65th & \parbox{0.31in}{\hfill \footnotesize 4.85} & 
50th & \parbox{0.31in}{\hfill \footnotesize 3.00} & 
35th & \parbox{0.31in}{\hfill \footnotesize 1.50} & 
20th & \parbox{0.31in}{\hfill \footnotesize 0.50} & 
5th & \parbox{0.31in}{\hfill \footnotesize 0}\\
\hline 
90th & \parbox{0.31in}{\hfill \footnotesize 13.33} & 
75th & \parbox{0.31in}{\hfill \footnotesize 6.50} & 
60th & \parbox{0.31in}{\hfill \footnotesize 4.00} & 
45th & \parbox{0.31in}{\hfill \footnotesize 2.50} & 
30th & \parbox{0.31in}{\hfill \footnotesize 1.00} & 
15th & \parbox{0.31in}{\hfill \footnotesize 0.33} & 
0th & \parbox{0.31in}{\hfill \footnotesize 0}\\
\hline 
\end{tabular}

\ 

\underline{\hspace*{4in}} 

\

Table 1.29: Career totals  for single-authored pages. 

\begin{tabular}{|r|c||r|c||r|c||r|c||r|c||r|c||r|c|}
\hline 
100th & \parbox{0.31in}{\hfill \footnotesize 780} & 
85th & \parbox{0.31in}{\hfill \footnotesize 64} & 
70th & \parbox{0.31in}{\hfill \footnotesize 35.5} & 
55th & \parbox{0.31in}{\hfill \footnotesize 21} & 
40th & \parbox{0.31in}{\hfill \footnotesize 11} & 
25th & \parbox{0.31in}{\hfill \footnotesize 0} & 
10th & \parbox{0.31in}{\hfill \footnotesize 0}\\
\hline 
95th & {\footnotesize 133.25} & 
80th & \parbox{0.31in}{\hfill \footnotesize 52} & 
65th & \parbox{0.31in}{\hfill \footnotesize 30.25} & 
50th & \parbox{0.31in}{\hfill \footnotesize 16.5} & 
35th & \parbox{0.31in}{\hfill \footnotesize 5} & 
20th & \parbox{0.31in}{\hfill \footnotesize 0} & 
5th & \parbox{0.31in}{\hfill \footnotesize 0}\\
\hline 
90th & \parbox{0.31in}{\hfill \footnotesize 88} & 
75th & \parbox{0.31in}{\hfill \footnotesize 42} & 
60th & \parbox{0.31in}{\hfill \footnotesize 25} & 
45th & \parbox{0.31in}{\hfill \footnotesize 14} & 
30th & \parbox{0.31in}{\hfill \footnotesize 0} & 
15th & \parbox{0.31in}{\hfill \footnotesize 0} & 
0th & \parbox{0.31in}{\hfill \footnotesize 0}\\
\hline 
\end{tabular}

\ 

Table 1.30: Career totals  for collaborative pages. 

\begin{tabular}{|r|c||r|c||r|c||r|c||r|c||r|c||r|c|}
\hline 
100th & \parbox{0.31in}{\hfill \footnotesize 2239} & 
85th & \parbox{0.31in}{\hfill \footnotesize 139.5} & 
70th & \parbox{0.31in}{\hfill \footnotesize 63.5} & 
55th & \parbox{0.31in}{\hfill \footnotesize 32} & 
40th & \parbox{0.31in}{\hfill \footnotesize 13} & 
25th & \parbox{0.31in}{\hfill \footnotesize 0} & 
10th & \parbox{0.31in}{\hfill \footnotesize 0}\\
\hline 
95th & \parbox{0.31in}{\hfill \footnotesize 306} & 
80th & \parbox{0.31in}{\hfill \footnotesize 112} & 
65th & \parbox{0.31in}{\hfill \footnotesize 50.25} & 
50th & \parbox{0.31in}{\hfill \footnotesize 25.5} & 
35th & \parbox{0.31in}{\hfill \footnotesize 8} & 
20th & \parbox{0.31in}{\hfill \footnotesize 0} & 
5th & \parbox{0.31in}{\hfill \footnotesize 0}\\
\hline 
90th & \parbox{0.31in}{\hfill \footnotesize 203.5} & 
75th & \parbox{0.31in}{\hfill \footnotesize 81.25} & 
60th & \parbox{0.31in}{\hfill \footnotesize 43} & 
45th & \parbox{0.31in}{\hfill \footnotesize 17} & 
30th & \parbox{0.31in}{\hfill \footnotesize 3} & 
15th & \parbox{0.31in}{\hfill \footnotesize 0} & 
0th & \parbox{0.31in}{\hfill \footnotesize 0}\\
\hline 
\end{tabular}

\ 

Table 1.31: Career totals  for attributable collaborative pages. 

\begin{tabular}{|r|c||r|c||r|c||r|c||r|c||r|c||r|c|}
\hline 
100th & {\footnotesize 890.10} & 
85th & \parbox{0.31in}{\hfill \footnotesize 60.92} & 
70th & \parbox{0.31in}{\hfill \footnotesize 29.00} & 
55th & \parbox{0.31in}{\hfill \footnotesize 13.63} & 
40th & \parbox{0.31in}{\hfill \footnotesize 6.00} & 
25th & \parbox{0.31in}{\hfill \footnotesize 0} & 
10th & \parbox{0.31in}{\hfill \footnotesize 0}\\
\hline 
95th & {\footnotesize 133.09} & 
80th & \parbox{0.31in}{\hfill \footnotesize 47.92} & 
65th & \parbox{0.31in}{\hfill \footnotesize 22.75} & 
50th & \parbox{0.31in}{\hfill \footnotesize 10.00} & 
35th & \parbox{0.31in}{\hfill \footnotesize 3.88} & 
20th & \parbox{0.31in}{\hfill \footnotesize 0} & 
5th & \parbox{0.31in}{\hfill \footnotesize 0}\\
\hline 
90th & \hspace*{0.01in} \parbox{0.31in}{\hfill \footnotesize 87.59} & 
75th & \parbox{0.31in}{\hfill \footnotesize 35.90} & 
60th & \parbox{0.31in}{\hfill \footnotesize 17.42} & 
45th & \parbox{0.31in}{\hfill \footnotesize 8.08} & 
30th & \parbox{0.31in}{\hfill \footnotesize 1.50} & 
15th & \parbox{0.31in}{\hfill \footnotesize 0} & 
0th & \parbox{0.31in}{\hfill \footnotesize 0}\\
\hline 
\end{tabular}

\ 

Table 1.32: Career totals  for attributable pages. 

\begin{tabular}{|r|c||r|c||r|c||r|c||r|c||r|c||r|c|}
\hline 
100th & {\scriptsize 1001.10} & 
85th & {\scriptsize 128.36} & 
70th & \parbox{0.31in}{\hfill \footnotesize 62.89} & 
55th & \parbox{0.31in}{\hfill \footnotesize 41.98} & 
40th & \parbox{0.31in}{\hfill \footnotesize 23.33} & 
25th & \parbox{0.31in}{\hfill \footnotesize 10.13} & 
10th & \parbox{0.31in}{\hfill \footnotesize 0}\\
\hline 
95th & {\scriptsize \hspace*{0.027in} 251.02} & 
80th & {\scriptsize 102.50} & 
65th & \parbox{0.31in}{\hfill \footnotesize 54.96} & 
50th & \parbox{0.31in}{\hfill \footnotesize 36.67} & 
35th & \parbox{0.31in}{\hfill \footnotesize 18.00} & 
20th & \parbox{0.31in}{\hfill \footnotesize 6.50} & 
5th & \parbox{0.31in}{\hfill \footnotesize 0}\\
\hline 
90th & {\scriptsize \hspace*{0.027in} 163.50} & 
75th & {\scriptsize \hspace*{0.032in} 75.00} & 
60th & \parbox{0.31in}{\hfill \footnotesize 48.00} & 
45th & \parbox{0.31in}{\hfill \footnotesize 31.00} & 
30th & \parbox{0.31in}{\hfill \footnotesize 14.33} & 
15th & \parbox{0.31in}{\hfill \footnotesize 3.44} & 
0th & \parbox{0.31in}{\hfill \footnotesize 0}\\
\hline 
\end{tabular}
\end{center} 
}

\newpage 
{\small 
\begin{center} 
\fbox{Table 2.1: time frame = 2002-2006}
 
\ClarifyingRemark 

\vspace*{0.1in}
Table 2.1: Biennial averages of single-authored and collaborative 
papers by university, for the years 2002-2006.\\
Reported for each metric are the mean, median, and mean of the middle 
50\% for the eligible faculty at each school. 

\begin{tabular}{|c||c||c||c|}
\hline 
\rule[-3mm]{0mm}{8mm}
University & Faculty & 
\parbox{1.5in}{\hspace*{0.07in} Single-authored Papers} & 
\parbox{1.5in}{\hspace*{0.13in} Collaborative Papers}\\
\hline 
\hline 
\rule[-1mm]{0mm}{4.5mm}Appalachian State (NC) & 10 & 
\statset{.32}{0}{0} & 
\statset{.20}{0}{.04}\\
\hline 
\rule[-1mm]{0mm}{4.5mm}Arkansas State & 6 & 
\statset{.53}{.2}{.33} & 
\statset{.07}{0}{0}\\
\hline 
\rule[-1mm]{0mm}{4.5mm}Central Connecticut State & 8 & 
\statset{.05}{0}{0} & 
\statset{.85}{.6}{.53}\\
\hline 
\rule[-1mm]{0mm}{4.5mm}Central Missouri State & 8 & 
\statset{.10}{0}{0} & 
\statset{.40}{0}{0}\\
\hline 
\rule[-1mm]{0mm}{4.5mm}Chicago State & 10 & 
\statset{0}{0}{0} & 
\statset{.04}{0}{0}\\
\hline 
\rule[-1mm]{0mm}{4.5mm}Eastern Illinois & 12 & 
\statset{.07}{0}{0} & 
\statset{.30}{0}{0}\\
\hline 
\rule[-1mm]{0mm}{4.5mm}Eastern Kentucky & 14 & 
\statset{.23}{0}{0} & 
\statset{.26}{0}{0}\\
\hline 
\end{tabular}

\ 
\end{center}
}

\newpage 
{\small 
\begin{center} 
\fbox{Table 2.2: time frame = 2002-2006}
 
\ClarifyingRemark
 
\vspace*{0.1in}
Table 2.2: Biennial averages of attributable collaborative and 
attributable 
papers by university, for 2002-2006.\\
Reported for each metric are the mean, median, and mean of the middle 
50\% for the eligible faculty at each school. 

\begin{tabular}{|c||c||c||c|}
\hline 
\rule[-4mm]{0mm}{10mm}
University & Faculty & 
\parbox{1.5in}{\hspace*{0.36in}Attributable\\ \hspace*{0.16in}Collaborative Papers} & 
\parbox{1.5in}{\hspace*{0.18in} Attributable Papers}\\
\hline 
\hline 
\rule[-1mm]{0mm}{4.5mm}Appalachian State (NC) & 10 & 
\statset{.08}{0}{.01} & 
\statset{.40}{0}{.09}\\
\hline 
\rule[-1mm]{0mm}{4.5mm}Arkansas State & 6 & 
\statset{.03}{0}{0} & 
\statset{.57}{.20}{.37}\\
\hline 
\rule[-1mm]{0mm}{4.5mm}Central Connecticut State & 8 & 
\statset{.37}{.27}{.22} & 
\statset{.42}{.43}{.35}\\
\hline 
\rule[-1mm]{0mm}{4.5mm}Central Missouri State & 8 & 
\statset{.18}{0}{0} & 
\statset{.28}{0}{0}\\
\hline 
\rule[-1mm]{0mm}{4.5mm}Chicago State & 10 & 
\statset{.01}{0}{0} & 
\statset{.01}{0}{0}\\
\hline 
\rule[-1mm]{0mm}{4.5mm}Eastern Illinois & 12 & 
\statset{.14}{0}{0} & 
\statset{.21}{0}{0}\\
\hline 
\rule[-1mm]{0mm}{4.5mm}Eastern Kentucky & 14 & 
\statset{.10}{0}{0} & 
\statset{.33}{0}{.07}\\
\hline 
\end{tabular}

\ 
\end{center}
}

\newpage 
{\small 
\begin{center} 
\fbox{Table 2.3: time frame = 2002-2006}

\ClarifyingRemark
 
\vspace*{0.1in}
Table 2.3: Biennial averages of single-authored and collaborative 
pages by university, for the years 2002-2006.\\
Reported for each metric are the mean, median, and mean of the middle 
50\% for the eligible faculty at each school. 

\begin{tabular}{|c||c||c||c|}
\hline 
\rule[-3mm]{0mm}{8mm}
University & Faculty & 
\parbox{1.5in}{\hspace*{0.09in} Single-authored Pages} & 
\parbox{1.5in}{\hspace*{0.15in} Collaborative Pages}\\
\hline 
\hline 
\rule[-1mm]{0mm}{4.5mm}Appalachian State (NC) & 10 & 
\statset{2.80}{0}{0} & 
\statset{2.44}{0}{.64}\\
\hline 
\rule[-1mm]{0mm}{4.5mm}Arkansas State & 6 & 
\statset{6.67}{2.6}{4.27} & 
\statset{1.27}{0}{0}\\
\hline 
\rule[-1mm]{0mm}{4.5mm}Central Connecticut State & 8 & 
\statset{.60}{0}{0} & 
\statset{11.70}{7.2}{6.27}\\
\hline 
\rule[-1mm]{0mm}{4.5mm}Central Missouri State & 8 & 
\statset{.95}{0}{0} & 
\statset{6.05}{0}{0}\\
\hline 
\rule[-1mm]{0mm}{4.5mm}Chicago State & 10 & 
\statset{0}{0}{0} & 
\statset{.52}{0}{0}\\
\hline 
\rule[-1mm]{0mm}{4.5mm}Eastern Illinois & 12 & 
\statset{.63}{0}{0} & 
\statset{5.37}{0}{0}\\
\hline 
\rule[-1mm]{0mm}{4.5mm}Eastern Kentucky & 14 & 
\statset{2.80}{0}{0} & 
\statset{2.77}{0}{0}\\
\hline 
\end{tabular}

\ 
\end{center}
}

\newpage 
{\small 
\begin{center} 
\fbox{Table 2.4: time frame = 2002-2006}

\ClarifyingRemark
 
\vspace*{0.1in}
Table 2.4: Biennial averages of attributable collaborative and 
attributable 
pages by university, for 2002-2006.\\
Reported for each metric are the mean, median, and mean of the middle 
50\% for the eligible faculty at each school. 

\begin{tabular}{|c||c||c||c|}
\hline 
\rule[-4mm]{0mm}{10mm}
University & Faculty & 
\parbox{1.5in}{\hspace*{0.36in}Attributable\\ \hspace*{0.18in}Collaborative 
Pages} & 
\parbox{1.5in}{\hspace*{0.2in} Attributable Pages}\\
\hline 
\hline 
\rule[-1mm]{0mm}{4.5mm}Appalachian State (NC) & 10 & 
\statset{.93}{0}{.16} & 
\statset{3.73}{0}{1.09}\\
\hline 
\rule[-1mm]{0mm}{4.5mm}Arkansas State & 6 & 
\statset{.63}{0}{0} & 
\statset{7.30}{2.60}{4.27}\\
\hline 
\rule[-1mm]{0mm}{4.5mm}Central Connecticut State & 8 & 
\statset{5.15}{3.07}{2.67} & 
\statset{5.75}{4.83}{3.95}\\
\hline 
\rule[-1mm]{0mm}{4.5mm}Central Missouri State & 8 & 
\statset{2.64}{0}{0} & 
\statset{3.59}{0}{0}\\
\hline 
\rule[-1mm]{0mm}{4.5mm}Chicago State & 10 & 
\statset{.17}{0}{0} & 
\statset{.17}{0}{0}\\
\hline 
\rule[-1mm]{0mm}{4.5mm}Eastern Illinois & 12 & 
\statset{2.35}{0}{0} & 
\statset{2.98}{0}{0}\\
\hline 
\rule[-1mm]{0mm}{4.5mm}Eastern Kentucky & 14 & 
\statset{1.15}{0}{0} & 
\statset{3.95}{0}{.46}\\
\hline 
\end{tabular}

\ 
\end{center}
}

\newpage 
{\small 
\begin{center} 
\fbox{Table 2.5: time frame = best five years}
 
\ClarifyingRemark
 
\vspace*{0.1in}
Table 2.5: Biennial averages of single-authored and collaborative 
papers by university, for the best five years.\\
Reported for each metric are the mean, median, and mean of the middle 
50\% for the eligible faculty at each school. 

\begin{tabular}{|c||c||c||c|}
\hline 
\rule[-3mm]{0mm}{8mm}
University & Faculty & 
\parbox{1.5in}{\hspace*{0.07in} Single-authored Papers} & 
\parbox{1.5in}{\hspace*{0.13in} Collaborative Papers}\\
\hline 
\hline 
\rule[-1mm]{0mm}{4.5mm}Appalachian State (NC) & 10 & 
\statset{.84}{.6}{.64} & 
\statset{.60}{.4}{.48}\\
\hline 
\rule[-1mm]{0mm}{4.5mm}Arkansas State & 6 & 
\statset{.73}{.4}{.53} & 
\statset{.27}{.4}{.33}\\
\hline 
\rule[-1mm]{0mm}{4.5mm}Central Connecticut State & 8 & 
\statset{.50}{.4}{.27} & 
\statset{1.30}{1.2}{1.07}\\
\hline 
\rule[-1mm]{0mm}{4.5mm}Central Missouri State & 8 & 
\statset{.60}{.4}{.27} & 
\statset{1.80}{.8}{.67}\\
\hline 
\rule[-1mm]{0mm}{4.5mm}Chicago State & 10 & 
\statset{.24}{0}{0} & 
\statset{.20}{0}{.04}\\
\hline 
\rule[-1mm]{0mm}{4.5mm}Eastern Illinois & 12 & 
\statset{1.03}{1.0}{.70} & 
\statset{1.03}{.8}{.60}\\
\hline 
\rule[-1mm]{0mm}{4.5mm}Eastern Kentucky & 14 & 
\statset{.49}{.4}{.40} & 
\statset{.43}{0}{.14}\\
\hline 
\end{tabular}

\ 
\end{center}
}

\newpage 
{\small 
\begin{center} 
\fbox{Table 2.6: time frame = best five years}

\ClarifyingRemark
 
\vspace*{0.1in}
Table 2.6: Biennial averages of attributable collaborative and 
attributable 
papers by university, for the best five years.\\
Reported for each metric are the mean, median, and mean of the middle 
50\% for the eligible faculty at each school. 

\begin{tabular}{|c||c||c||c|}
\hline 
\rule[-4mm]{0mm}{10mm}
University & Faculty & 
\parbox{1.5in}{\hspace*{0.36in}Attributable\\ \hspace*{0.16in}Collaborative Papers} & 
\parbox{1.5in}{\hspace*{0.18in} Attributable Papers}\\
\hline 
\hline 
\rule[-1mm]{0mm}{4.5mm}Appalachian State (NC) & 10 & 
\statset{.24}{.17}{.19} & 
\statset{.94}{.77}{.77}\\
\hline 
\rule[-1mm]{0mm}{4.5mm}Arkansas State & 6 & 
\statset{.13}{.20}{.17} & 
\statset{.83}{.60}{.70}\\
\hline 
\rule[-1mm]{0mm}{4.5mm}Central Connecticut State & 8 & 
\statset{.61}{.60}{.51} & 
\statset{.91}{.77}{.69}\\
\hline 
\rule[-1mm]{0mm}{4.5mm}Central Missouri State & 8 & 
\statset{.79}{.37}{.31} & 
\statset{1.11}{.80}{.60}\\
\hline 
\rule[-1mm]{0mm}{4.5mm}Chicago State & 10 & 
\statset{.09}{0}{.01} & 
\statset{.29}{0}{.04}\\
\hline 
\rule[-1mm]{0mm}{4.5mm}Eastern Illinois & 12 & 
\statset{.47}{.30}{.25} & 
\statset{1.23}{1.00}{.80}\\
\hline 
\rule[-1mm]{0mm}{4.5mm}Eastern Kentucky & 14 & 
\statset{.19}{0}{.05} & 
\statset{.61}{.40}{.49}\\
\hline 
\end{tabular}

\ 
\end{center}
}

\newpage 
{\small 
\begin{center} 
\fbox{Table 2.7: time frame = best five years}

\ClarifyingRemark
 
\vspace*{0.1in}
Table 2.7: Biennial averages of single-authored and collaborative 
pages by university, for the best five years.\\
Reported for each metric are the mean, median, and mean of the middle 
50\% for the eligible faculty at each school. 

\begin{tabular}{|c||c||c||c|}
\hline 
\rule[-3mm]{0mm}{8mm}
University & Faculty & 
\parbox{1.5in}{\hspace*{0.09in} Single-authored Pages} & 
\parbox{1.5in}{\hspace*{0.15in} Collaborative Pages}\\
\hline 
\hline 
\rule[-1mm]{0mm}{4.5mm}Appalachian State (NC) & 10 & 
\statset{10.00}{10.4}{9.28} & 
\statset{10.48}{10.8}{10.96}\\
\hline 
\rule[-1mm]{0mm}{4.5mm}Arkansas State & 6 & 
\statset{14.07}{10.6}{11.27} & 
\statset{6.67}{6.2}{5.40}\\
\hline 
\rule[-1mm]{0mm}{4.5mm}Central Connecticut State & 8 & 
\statset{4.85}{3.2}{2.13} & 
\statset{17.55}{15.2}{13.07}\\
\hline 
\rule[-1mm]{0mm}{4.5mm}Central Missouri State & 8 & 
\statset{4.75}{2.0}{1.33} & 
\statset{14.35}{7.6}{5.87}\\
\hline 
\rule[-1mm]{0mm}{4.5mm}Chicago State & 10 & 
\statset{2.48}{0}{0} & 
\statset{3.68}{0}{.52}\\
\hline 
\rule[-1mm]{0mm}{4.5mm}Eastern Illinois & 12 & 
\statset{10.07}{8.8}{7.80} & 
\statset{14.73}{8.6}{7.60}\\
\hline 
\rule[-1mm]{0mm}{4.5mm}Eastern Kentucky & 14 & 
\statset{5.54}{4.6}{4.17} & 
\statset{5.49}{0}{1.37}\\
\hline 
\end{tabular}

\ 
\end{center}
}

\newpage 
{\small 
\begin{center} 
\fbox{Table 2.8: time frame = best five years}
 
\ClarifyingRemark
 
\vspace*{0.1in}
Table 2.8: Biennial averages of attributable collaborative and 
attributable 
pages by university, for the best five years.\\
Reported for each metric are the mean, median, and mean of the middle 
50\% for the eligible faculty at each school. 

\begin{tabular}{|c||c||c||c|}
\hline 
\rule[-4mm]{0mm}{10mm}
University & Faculty & 
\parbox{1.5in}{\hspace*{0.36in}Attributable\\ \hspace*{0.18in}Collaborative 
Pages} & 
\parbox{1.5in}{\hspace*{0.2in} Attributable Pages}\\
\hline 
\hline 
\rule[-1mm]{0mm}{4.5mm}Appalachian State (NC) & 10 & 
\statset{4.11}{3.93}{4.14} & 
\statset{13.55}{12.40}{12.59}\\
\hline 
\rule[-1mm]{0mm}{4.5mm}Arkansas State & 6 & 
\statset{3.33}{3.10}{2.70} & 
\statset{15.73}{13.70}{13.97}\\
\hline 
\rule[-1mm]{0mm}{4.5mm}Central Connecticut State & 8 & 
\statset{8.13}{7.20}{6.11} & 
\statset{11.26}{10.40}{8.84}\\
\hline 
\rule[-1mm]{0mm}{4.5mm}Central Missouri State & 8 & 
\statset{6.04}{3.07}{2.44} & 
\statset{9.26}{6.40}{4.98}\\
\hline 
\rule[-1mm]{0mm}{4.5mm}Chicago State & 10 & 
\statset{1.23}{0}{.17} & 
\statset{3.71}{0}{.92}\\
\hline 
\rule[-1mm]{0mm}{4.5mm}Eastern Illinois & 12 & 
\statset{6.42}{4.30}{3.63} & 
\statset{12.72}{10.50}{9.35}\\
\hline 
\rule[-1mm]{0mm}{4.5mm}Eastern Kentucky & 14 & 
\statset{2.35}{0}{.53} & 
\statset{7.24}{6.20}{6.30}\\
\hline 
\end{tabular}

\ 
\end{center}
}

\newpage 
{\small 
\begin{center} 
\fbox{Table 2.9: time frame = best up-to-ten years}
 
\ClarifyingRemark
 
\vspace*{0.1in}
Table 2.9: Biennial averages of single-authored and collaborative 
papers by university, for the best up-to-ten years.\\
Reported for each metric are the mean, median, and mean of the middle 
50\% for the eligible faculty at each school. 

\begin{tabular}{|c||c||c||c|}
\hline 
\rule[-3mm]{0mm}{8mm}
University & Faculty & 
\parbox{1.5in}{\hspace*{0.07in} Single-authored Papers} & 
\parbox{1.5in}{\hspace*{0.13in} Collaborative Papers}\\
\hline 
\hline 
\rule[-1mm]{0mm}{4.5mm}Appalachian State (NC) & 10 & 
\statset{.59}{.33}{.38} & 
\statset{.35}{.23}{.25}\\
\hline 
\rule[-1mm]{0mm}{4.5mm}Arkansas State & 6 & 
\statset{.54}{.22}{.31} & 
\statset{.15}{.20}{.17}\\
\hline 
\rule[-1mm]{0mm}{4.5mm}Central Connecticut State & 8 & 
\statset{.32}{.20}{.13} & 
\statset{.96}{.78}{.65}\\
\hline 
\rule[-1mm]{0mm}{4.5mm}Central Missouri State & 8 & 
\statset{.33}{.20}{.13} & 
\statset{1.23}{.40}{.33}\\
\hline 
\rule[-1mm]{0mm}{4.5mm}Chicago State & 10 & 
\statset{.14}{0}{0} & 
\statset{.14}{0}{.02}\\
\hline 
\rule[-1mm]{0mm}{4.5mm}Eastern Illinois & 12 & 
\statset{.65}{.50}{.40} & 
\statset{.67}{.40}{.30}\\
\hline 
\rule[-1mm]{0mm}{4.5mm}Eastern Kentucky & 14 & 
\statset{.33}{.20}{.20} & 
\statset{.30}{0}{.07}\\
\hline 
\end{tabular}

\ 
\end{center}
}

\newpage 
{\small 
\begin{center} 
\fbox{Table 2.10: time frame = best up-to-ten years}

\ClarifyingRemark
 
\vspace*{0.1in}
Table 2.10: Biennial averages of attrib.\ collaborative and attributable 
papers by university, for the best up-to-ten years.\\
Reported for each metric are the mean, median, and mean of the middle 
50\% for the eligible faculty at each school. 

\begin{tabular}{|c||c||c||c|}
\hline 
\rule[-4mm]{0mm}{10mm}
University & Faculty & 
\parbox{1.5in}{\hspace*{0.36in}Attributable\\ \hspace*{0.16in}Collaborative Papers} & 
\parbox{1.5in}{\hspace*{0.18in} Attributable Papers}\\
\hline 
\hline 
\rule[-1mm]{0mm}{4.5mm}Appalachian State (NC) & 10 & 
\statset{.14}{.09}{.10} & 
\statset{.69}{.42}{.46}\\
\hline 
\rule[-1mm]{0mm}{4.5mm}Arkansas State & 6 & 
\statset{.07}{.10}{.08} & 
\statset{.60}{.33}{.41}\\
\hline 
\rule[-1mm]{0mm}{4.5mm}Central Connecticut State & 8 & 
\statset{.44}{.39}{.32} & 
\statset{.64}{.52}{.49}\\
\hline 
\rule[-1mm]{0mm}{4.5mm}Central Missouri State & 8 & 
\statset{.54}{.19}{.16} & 
\statset{.80}{.49}{.36}\\
\hline 
\rule[-1mm]{0mm}{4.5mm}Chicago State & 10 & 
\statset{.06}{0}{.01} & 
\statset{.18}{0}{.02}\\
\hline 
\rule[-1mm]{0mm}{4.5mm}Eastern Illinois & 12 & 
\statset{.30}{.15}{.13} & 
\statset{.84}{.55}{.46}\\
\hline 
\rule[-1mm]{0mm}{4.5mm}Eastern Kentucky & 14 & 
\statset{.13}{0}{.03} & 
\statset{.44}{.20}{.27}\\
\hline 
\end{tabular}

\ 
\end{center}
}

\newpage 
{\small 
\begin{center} 
\fbox{Table 2.11: time frame = best up-to-ten years}

\ClarifyingRemark
 
\vspace*{0.1in}
Table 2.11: Biennial averages of single-authored and collaborative 
pages by university, for the best up-to-ten years.\\
Reported for each metric are the mean, median, and mean of the middle 
50\% for the eligible faculty at each school. 

\begin{tabular}{|c||c||c||c|}
\hline 
\rule[-3mm]{0mm}{8mm}
University & Faculty & 
\parbox{1.5in}{\hspace*{0.09in} Single-authored Pages} & 
\parbox{1.5in}{\hspace*{0.15in} Collaborative Pages}\\
\hline 
\hline 
\rule[-1mm]{0mm}{4.5mm}Appalachian State (NC) & 10 & 
\statset{6.21}{5.20}{5.36} & 
\statset{5.99}{6.03}{5.85}\\
\hline 
\rule[-1mm]{0mm}{4.5mm}Arkansas State & 6 & 
\statset{8.86}{5.89}{6.55} & 
\statset{3.56}{3.30}{2.99}\\
\hline 
\rule[-1mm]{0mm}{4.5mm}Central Connecticut State & 8 & 
\statset{3.19}{2.50}{1.67} & 
\statset{12.39}{10.00}{8.13}\\
\hline 
\rule[-1mm]{0mm}{4.5mm}Central Missouri State & 8 & 
\statset{2.78}{1.00}{.67} & 
\statset{10.25}{3.80}{2.93}\\
\hline 
\rule[-1mm]{0mm}{4.5mm}Chicago State & 10 & 
\statset{1.50}{0}{0} & 
\statset{2.32}{0}{.26}\\
\hline 
\rule[-1mm]{0mm}{4.5mm}Eastern Illinois & 12 & 
\statset{6.13}{5.00}{4.20} & 
\statset{10.22}{4.30}{3.80}\\
\hline 
\rule[-1mm]{0mm}{4.5mm}Eastern Kentucky & 14 & 
\statset{4.00}{2.30}{2.09} & 
\statset{3.31}{0}{.86}\\
\hline 
\end{tabular}

\ 
\end{center}
}

\newpage 
{\small 
\begin{center} 
\fbox{Table 2.12: time frame = best up-to-ten years}
 
\ClarifyingRemark
 
\vspace*{0.1in}
Table 2.12: Biennial averages of attrib.\ collaborative and attributable 
pages by university, for the best up-to-ten years.\\
Reported for each metric are the mean, median, and mean of the middle 
50\% for the eligible faculty at each school. 

%

\ 
\end{center}
}

\newpage 
{\small 
\begin{center} 
\fbox{Table 2.13: time frame = career}
 
\ClarifyingRemark

\vspace*{0.1in}
Table 2.13: Career totals for single-authored and collaborative 
papers by university.\\
Reported for each metric are the mean, median, and mean of the middle 
50\% for the eligible faculty at each school. 

\begin{tabular}{|c||c||c||c|}
\hline 
\rule[-3mm]{0mm}{8mm}
University & Faculty & 
\parbox{1.5in}{\hspace*{0.07in} Single-authored Papers} & 
\parbox{1.5in}{\hspace*{0.13in} Collaborative Papers}\\
\hline 
\hline 
\rule[-1mm]{0mm}{4.5mm}Appalachian State (NC) & 10 & 
\statset{3.50}{1.5}{1.60} & 
\statset{2.00}{1.0}{1.20}\\
\hline 
\rule[-1mm]{0mm}{4.5mm}Arkansas State & 6 & 
\statset{2.67}{1.0}{1.33} & 
\statset{.67}{1.0}{.83}\\
\hline 
\rule[-1mm]{0mm}{4.5mm}Central Connecticut State & 8 & 
\statset{1.38}{1.0}{.67} & 
\statset{5.13}{4.0}{3.67}\\
\hline 
\rule[-1mm]{0mm}{4.5mm}Central Missouri State & 8 & 
\statset{1.75}{1.0}{.67} & 
\statset{7.88}{2.0}{1.67}\\
\hline 
\rule[-1mm]{0mm}{4.5mm}Chicago State & 10 & 
\statset{.70}{0}{0} & 
\statset{.70}{0}{.10}\\
\hline 
\rule[-1mm]{0mm}{4.5mm}Eastern Illinois & 12 & 
\statset{3.75}{3.0}{2.25} & 
\statset{4.83}{2.0}{1.75}\\
\hline 
\rule[-1mm]{0mm}{4.5mm}Eastern Kentucky & 14 & 
\statset{1.29}{1.0}{1.00} & 
\statset{1.64}{0}{.36}\\
\hline 
\end{tabular}

\ 
\end{center}
}

\newpage 
{\small 
\begin{center} 
\fbox{Table 2.14: time frame = career}
 
\ClarifyingRemark
 
\vspace*{0.1in}
Table 2.14: Career totals for attributable collaborative and 
attributable 
papers by university.\\
Reported for each metric are the mean, median, and mean of the middle 
50\% for the eligible faculty at each school. 

%

\ 
\end{center}
}

\newpage 
{\small 
\begin{center} 
\fbox{Table 2.15: time frame = career}

\ClarifyingRemark

\vspace*{0.1in}
Table 2.15: Career totals for single-authored and collaborative 
pages by university.\\
Reported for each metric are the mean, median, and mean of the middle 
50\% for the eligible faculty at each school. 

%

\ 
\end{center}
}

\newpage 
{\small 
\begin{center} 
\fbox{Table 2.16: time frame = career}
 
\ClarifyingRemark

\vspace*{0.1in}
Table 2.16: Career totals for attributable collaborative and 
attributable 
pages by university.\\
Reported for each metric are the mean, median, and mean of the middle 
50\% for the eligible faculty at each school. 

%

\ 
\end{center}
}

\newpage 
\def\refname{\normalsize References.}
\renewcommand{\baselinestretch}{1}

\end{document}